\documentclass[11pt,a4paper]{amsart}
\usepackage{amscd,amssymb,graphics,color}

\usepackage{mathrsfs}
\newtheorem{theorem}{Theorem}[section]
\newtheorem{lemma}[theorem]{Lemma}
\newtheorem{proposition}[theorem]{Proposition}
\newtheorem{corollary}[theorem]{Corollary}

\theoremstyle{definition}
\newtheorem{definition}[theorem]{Definition}
\newtheorem{construction}[theorem]{Construction}
\newtheorem{example}[theorem]{Example}

\theoremstyle{remark}
\newtheorem{remark}[theorem]{Remark}

\newcommand{\NN} {\mathbb{N}}
\newcommand{\ZZ} {\mathbb{Z}}
\newcommand{\QQ} {\mathbb{Q}}
\newcommand{\RR} {\mathbb{R}}
\newcommand{\CC} {\mathbb{C}}
\newcommand{\FF} {\mathbb{F}}
\newcommand{\PP} {\mathbb{P}}
\renewcommand{\AA} {\mathbb{A}}
\newcommand{\GG} {\mathbb{G}}

\newcommand {\shE} {\mathcal{E}}
\newcommand {\shF} {\mathcal{F}}

\newcommand {\shM} {\mathcal{M}}

\newcommand {\shN} {\mathcal{N}}

\newcommand {\A} {\mathfrak{A}}

\newcommand {\Aut} {\operatorname{Aut}}

\newcommand {\codim} {\operatorname{codim}}

\newcommand {\coker} {\operatorname{coker}}

\newcommand {\dlog} {\operatorname{dlog}}

\newcommand {\eps} {\varepsilon}

\newcommand {\Ext} {\operatorname{Ext}}

\newcommand {\Hom} {\operatorname{Hom}}

\newcommand {\im} {\operatorname{im}}

\newcommand {\Int} {\operatorname{int}}

\newcommand {\kk} {\Bbbk}
\renewcommand {\ker } {\operatorname{ker}}

\newcommand {\LL} {\mathfrak{L}}

\newcommand {\lra} {\longrightarrow}

\newcommand {\M} {\mathcal{M}}
\newcommand {\map} {\operatorname{Map}}

\newcommand {\maxid} {\mathfrak{m}}

\newcommand {\mult} {\operatorname{mult}}

\renewcommand{\O} {\mathcal{O}}

\renewcommand {\P} {\mathscr{P}}

\newcommand {\sing} {\mathrm{sing}}

\newcommand {\Spec} {\operatorname{Spec}}

\newcommand {\T} {\mathfrak{T}}

\newcommand{\bbfamily}{\fontencoding{U}\fontfamily{bbold}\selectfont}
\newcommand{\textbb}[1]{{\bbfamily#1}}
\newcommand {\lfor} {\mbox{\rm\textbb{[}}}
\newcommand {\rfor} {\mbox{\rm\textbb{]}}}

\def\mydate{\ifcase\month \or January\or February\or March\or
April\or May\or June\or July\or August\or September\or October\or 
November\or December\fi \space\number\day,\space\number\year}

\begin{document}

\title[Toric varieties and tropical curves]{Toric degenerations of
toric varieties and tropical curves}

\author{Takeo Nishinou} \address{Department of Mathematics, Faculty of
Science, Kyoto University, Kitashirakawa, Kyoto, 606-8224, Japan}
\curraddr{Mathematisches Institut,
Albert-Ludwigs-Universit\"at, Eckerstra\ss e~1, 79104 Freiburg,
Germany}
\email{nishinou@kusm.kyoto-u.ac.jp}
\author{Bernd Siebert} \address{Mathematisches Institut,
Albert-Ludwigs-Universit\"at, Eckerstra\ss e~1, 79104 Freiburg,
Germany}
\email{bernd.siebert@math.uni-freiburg.de}

\date{December 15, 2005}

\begin{abstract}
We show that the counting of rational curves on a complete toric
variety that are in general position to the toric prime divisors
coincides with the counting of certain tropical curves. The proof is
algebraic-geometric and relies on degeneration techniques and log
deformation theory.
\end{abstract}

\maketitle
\tableofcontents

\section*{Introduction.}

The advent of mirror symmetry has rekindled the interest in the
classical enumerative question of counting algebraic curves on a
smooth projective variety. In fact, the unexpected relation of this
problem with period computations for the mirror family \cite{candelas
etal} really caused the great mathematical interest in mirror
symmetry. While meanwhile verified mathematically by direct
computation \cite{givental} (see also \cite{lianliuyau}) no geometric
explanation for the validity of these computations have been found
yet.

In connection with the mirror symmetry program with M.~Gross
\cite{announce}, \cite{logmirror}, the second author has long
suspected that this counting problem should have a purely
combinatorial counterpart in integral affine geometry (cf.\ also
\cite{fukaya oh}, \cite{kontsoib}, and the work on Gromov-Witten
invariants for degenerations \cite{li ruan}, \cite{ionel parker},
\cite{junli2}, \cite{logGW}). Moreover, the Legendre dual
interpretation of the combinatorial counting problem should be related
to period computations on the mirror, thus eventually giving a
satisfactory explanation for the mirror computation in \cite{candelas
etal}.

Recently G.~Mikhalkin used certain piecewise linear, one-dimensional
objects in $\QQ^d$, called tropical curves, to give a count of
algebraic curves on toric surfaces \cite{mikhalkin}. As tropical
curves make perfectly sense in the integral affine context, it is
then natural to look at the counting of tropical curves from the
degeneration point of view of \cite{logmirror}. The pleasant surprise
was that not only do tropical curves come up very naturally as dual
intersection graphs of certain transverse curves on the central
fiber, but also that our logarithmic techniques give more transparent
and robust proofs in the toric case. In particular, as logarithmic
deformation theory replaces the hypersurface method of patchworking
used in \cite{mikhalkin}, our results are not limited to two
dimensions. We also believe that this technical robustness will
be essential in generalizations to the Calabi-Yau case. The purpose
of this paper is thus not only a generalization of the results of
\cite{mikhalkin} to arbitrary dimensions, but also to lay the
foundations for a treatment of the Calabi-Yau case in the framework of
\cite{logmirror}.

The main result is Theorem~\ref{main theorem}. It says how exactly
tropical curves count algebraic curves on an arbitrary complete toric
variety. The formulation of this theorem involves some concepts
introduced in the first three sections. We therefore devote only
Section~8 to statement and discussion of this theorem. This section
also ties together the various strands of proof given in the other
sections.

During the preparation of this paper we learned that G.~Mikhalkin has
also announced a generalization of his work to higher dimensions,
using symplectic techniques. We apologize for any duplication of
results. We also thank him for pointing out an elucidating example
concerning the multiplicity definition for tropical curves in higher
dimensions, which had not been worked out properly by us at that time.
\smallskip

\noindent
\emph{Conventions.}\ \ 
We work in the category of schemes of finite type over an
algebraically closed field $\kk$ of characteristic $0$. Throughout the
paper $N$ is a free abelian group of rank $n\ge 2$ and $N_\QQ=
N\otimes_\ZZ \QQ$, $N_\RR=N\otimes_\ZZ \RR$. For toric geometry fans
will be defined in $N_\QQ$ or in $N_\QQ\times\QQ$. If $\Sigma$ is a
fan in $N$ then $X(\Sigma)$ is the associated toric $\kk$-variety with
big torus $\Int X(\Sigma)\simeq \GG(N)\subset X(\Sigma)$. Its
complement is sometimes referred to as the \emph{toric boundary} of
$X(\Sigma)$. For a polyhedral cone $\sigma\subset N_\QQ$ the union of
proper faces of $\sigma$ is denoted $\partial\sigma$, and $\Int
\sigma= \sigma\setminus \partial\sigma$ is the relative interior.
Dually we have $M=\Hom(N,\ZZ)$ and $M_\QQ=M\otimes_\ZZ\QQ$, and the
dual pairing $M\times N\to \QQ$ is denoted by brackets $\langle\
\;,\;\rangle$. The notation for the algebraic torus $\Spec \kk[M]$ is
$\GG(N)$. Then $N$ and $M$ can be identified with the space of
one-parameter subgroups $\GG_m\to \GG(N)$ and of characters
$\GG(N)\to\GG_m$, respectively. This motivates the notation $\chi^m$
for the monomial in $\kk[M]$, or in any subring, corresponding to
$m\in M$. If $\Xi\subset N_\QQ$ is a subset $L(\Xi)\subset N_\QQ$
denotes the linear subspace spanned by differences $v-w$ for $v,w\in
\Xi$, and $C(\Xi)\subset N_\QQ\times\QQ$ is the closure of the convex
hull of $\QQ_{\ge 0}\cdot (A\times\{1\})$. In particular, if $A$ is an
affine subspace then $L(A)\subset N_\QQ$ is the associated linear
space and $LC(A):=L(C(A))\subset \NN_\QQ\times\QQ$ is the linear
closure of $A\times\{1\}$. The natural numbers $\NN$ include $0$.

\section{Tropical curves}

Recall the definition of tropical curve given in \cite{mikhalkin},
Definition~2.2. Let $\overline{\Gamma}$ be (the geometric realization
of) a weighted, connected finite graph \emph{without divalent
vertices}. Its sets of vertices and edges are denoted
$\overline{\Gamma}^{[0]}$ and $\overline{\Gamma}^{[1]}$ respectively,
and $w_{\overline{\Gamma}}: \overline{\Gamma}^{[1]}\to \NN\setminus
\{0\}$ is the weight function.  An edge $E\in\overline{\Gamma}^{[1]}$
has adjacent vertices $\partial E=\{V_1,V_2\}$. Let
$\overline{\Gamma}^{[0]}_\infty \subset \overline{\Gamma}^{[0]}$ be
the set of one-valent vertices. We set
\[
\Gamma = \overline{\Gamma} \setminus \overline{\Gamma}^{[0]}_\infty.
\]
By abuse of notation let us refer to such an object as a (weighted)
\emph{open graph} with vertices, edges and weight function
$\Gamma^{[0]}$, $\Gamma^{[1]}$, $w_{\Gamma}$. Now some edges may be
non-compact, and these are called \emph{unbounded edges}. Write
$\Gamma_\infty^{[1]} \subset \Gamma^{[1]}$ for the subset of
unbounded edges. The product of all weights of bounded edges is the
\emph{total inner weight} of $\Gamma$
\[
w(\Gamma):= \prod_{E\in \Gamma^{[1]}\setminus \Gamma^{[1]}_\infty}
w_\Gamma(E).
\]
The set of \emph{flags} of $\Gamma$ is
$F\Gamma=\{(V,E)\,|\, V\in\partial E\}$. Because every unbounded edge
has only one adjacent vertex we sometimes consider
$\Gamma_\infty^{[1]}$ as subset of $F\Gamma$.

\begin{definition}
A \emph{parameterized tropical curve} in $N_\QQ$ is a proper map $h:
\Gamma \to N_\RR= N_\QQ\otimes_\QQ\RR$ satisfying the following
conditions.
\begin{enumerate}
\item[(i)] For every edge $E \subset \Gamma$ the restriction $h|_E$ is
an embedding with image $h(E)$ contained in an affine line with
rational slope.
\item[(ii)] For every vertex $V \in \Gamma$ we have $h(V)\in N_\QQ$
and the following \emph{balancing condition} holds. Let $E_1, \dots,
E_m \in \Gamma^{[1]}$ be the edges adjacent to $V$, and let $u_i \in
N$ be the primitive integral vector emanating from $h(V)$ in the
direction of $h(E_i)$. Then
\begin{eqnarray}\label{balancing condition}
\sum_{j=1}^m w(E_j) u_j = 0.
\end{eqnarray}
\end{enumerate}
An \emph{isomorphism} of tropical curves $h: \Gamma \to N_\RR$ and
$h': \Gamma' \to N_\RR$ is a homeomorphism $\Phi: \Gamma \to \Gamma'$
respecting the weights of the edges and such that $h = h' \circ
\Phi$. A \emph{tropical curve} is an isomorphism class of
parameterized tropical curves. Note that the isomorphism class of a
tropical curve is determined by the weighted graph $\Gamma$, a map
$\Gamma^{[0]}\to N_\QQ$ telling the images of the vertices, and a map
$\Gamma^{[1]}_\infty \to N_\QQ$ for the slopes of the unbounded
edges. The \emph{genus} of a tropical curve $h:\Gamma\to N_\RR$ is
the first Betti number of $\Gamma$. A \emph{rational} tropical curve
is a tropical curve of genus $0$.
\qed
\end{definition}

\begin{remark}
From a systematic point of view one might also want to allow
contracted edges, as in the latest version of \cite{mikhalkin}. In
strict analogy with stable maps one should then also introduce a
stability condition involving marked edges. For the unobstructed
cases treated in this paper contracted edges are irrelevant.
\end{remark}

Implicit in (ii) is the definition of a map $u:F\Gamma\to N$ sending
a flag $(V,E)$ to the primitive integral vector $u_{(V,E)}\in N$
emanating from $V$ in the direction of $h(E)$. Note that
$u_{(V_1,E)}=-u_{(V_2,E)}$ if $\partial E=\{V_1,V_2\}$. With this
notation the balancing condition at $V\in\Gamma^{[0]}$ reads
$\sum_{\{E\in\Gamma^{[1]}\,|\, V\in\partial E\}} w(E)\cdot
u_{(V,E)}=0$.

To impose incidence conditions on our tropical curve we need to
distinguish some edges. Let $l\in \NN$. An \emph{$l$-marked} tropical
curve is a tropical curve $h:\Gamma\to N_\RR$ together with a choice
of $l$ edges $\mathbf{E}=(E_1,\ldots, E_l)\subset (\Gamma^{[1]})^l$.
In this definition we do not assume the $E_i$ to be pairwise
distinct. The notation is $(\Gamma,\mathbf{E}, h)$. The \emph{total
marked weight} of the marked weighted graph $(\Gamma,\mathbf{E})$ is
\[
w(\Gamma,\mathbf{E}):= w(\Gamma)\cdot\prod_{i=1}^l w(E_i).
\]
This is the number that will later be relevant for enumerative
geometry, see Proposition~\ref{log structures on degenerate curves}
and Section~\ref{section main theorem}. The \emph{type} of
$(\Gamma,\mathbf{E},h)$ is the marked graph $(\Gamma,\mathbf{E})$
together with the map $u:F\Gamma \to N$. For given type
$(\Gamma,\mathbf{E},u)$ denote by $\T_{(\Gamma,\mathbf{E},u)}
\simeq\T_{(\Gamma,u)}$ the space of isomorphism classes of (marked)
tropical curves of this type. If non-empty it is a manifold of
dimension
\[
\dim \T_{(\Gamma,u)}\ge e+(n-3)(1-g)-\operatorname{ov}(\Gamma),
\]
where $e$ is the number of unbounded edges, $g$ is the genus of
$\Gamma$ and $\operatorname{ov}(\Gamma)$ is the \emph{overvalence} of
$\Gamma$ (\cite{mikhalkin}, Proposition~2.14). The overvalence
measures the difference from $\Gamma$ to a trivalent graph:
\[
\operatorname{ov}(\Gamma) =\sum_{V\in\Gamma^{[0]}}
\Big(\sharp \big\{E\in\Gamma^{[1]} \,\big|\, (V,E)\in F\Gamma\big\}
-3\Big).
\]
Local coordinates are given by the position of one vertex
and the lengths of the edges of a spanning tree, subject to
linear relations coming from the remaining edges. In particular, for $g=0$
\begin{eqnarray}\label{dim space of tropical curves}
\dim \T_{(\Gamma,u)}= e+n-3-\operatorname{ov}(\Gamma)
\le e+n-3,
\end{eqnarray}
with equality iff $\Gamma$ is trivalent, and $e-3- \operatorname{ov}
(\Gamma)$ equals the number of bounded edges of $\Gamma$.

\begin{definition}
For $\mathbf{d}=(d_1,\ldots,d_l)\in\NN^l$ an \emph{affine
constraint} of codimension $\mathbf{d}$ is an $l$-tuple
$\mathbf{A}=(A_1,\ldots, A_l)$ of affine subspaces $A_i\subset N_\QQ$
with $\dim A_i= n-d_i-1$. An $l$-marked tropical curve
$(\Gamma,\mathbf{E}, h)$ \emph{matches} the affine constraint
$\mathbf{A}$ if
\[
h(E_i)\cap A_i\neq \emptyset,\quad i=1,\ldots,l.
\]
The space of $l$-marked tropical curves of type
$(\Gamma,\mathbf{E},u)$ matching $\mathbf{A}$ is denoted
$\T_{(\Gamma,\mathbf{E},u)}(\mathbf{A})$.
\end{definition}

The \emph{degree} of a type $(\Gamma,\mathbf{E},u)$ of tropical
curves is the function $N\setminus \{0\}\to \NN$ with finite support
defined by
\[
\Delta(\Gamma,\mathbf{E},u)(v) =\Delta(\Gamma,u)(v)
:=\sharp\big\{ (V,E)\in\Gamma_\infty^{[1]}
\,\big|\, w(E)\cdot u_{(V,E)}=v\big\},
\]
where we consider $\Gamma_\infty^{[1]}$ as a subset of $F\Gamma$. In
other words, it is the abstract set of directions of unbounded edges
together with their weights, with repetitions allowed. The degree of
a (marked) tropical curve is the degree of its type. For $g\in \NN$
and $\Delta\in\map(N\setminus\{0\},\NN)$ the set of
$l$-marked tropical curves of genus $g$ and degree $\Delta$ is denoted
$\T_{g,l,\Delta}$. For the subset matching an affine constraint
$\mathbf A$ the notation is $\T_{g,l,\Delta}(\mathbf{A})$.
For $\Delta\in\map(N\setminus \{0\},\NN)$ with finite support define
$|\Delta|=\sum_{v\in N\setminus\{0\}} \Delta(v)$. A tropical curve of
degree $\Delta$ has $|\Delta|$ unbounded edges (unweighted count).

Note that by the balancing condition~(\ref{balancing condition})
there are no tropical curves of degree $\Delta\in \map(N\setminus
\{0\},\NN)$ unless $\sum_{v\in N\setminus \{0\}} \Delta(v)\cdot v=0$. The
enumerative count on the affine side is of tropical curves of fixed
degree and matching a given general affine constraint. To make this
precise the next section discusses transversality in this context.
Subsequent sections treat the toric side, with the degree selecting
the homology class of curves to be considered, and the affine
constraints corresponding to orbits of subgroups to be intersected by
the algebraic curves.

\section{Affine transversality}

\begin{proposition}\label{finiteness of types}
For any $\Delta\in \map(N\setminus\{0\},\NN)$ and any $g\in\NN$ there
are only finitely many types of tropical curves of degree $\Delta$
and genus $g$, that is, the set
\[
\big\{(\Gamma,u)\text{ type of tropical curve}\,\big|\,
\Delta(\Gamma,u)=\Delta,\ g(\Gamma,u)=g,\,
\T_{(\Gamma,u)}\neq\emptyset\big\}
\]
is finite.
\end{proposition}
\proof
Let $(\Gamma,u)$ have degree $\Delta$ and genus $g$. Then the number
of unbounded edges of $\Gamma$ is equal to $|\Delta|$. Hence $\Delta$
must have finite support. Moreover, for each $e$ there are only
finitely many graphs $\overline\Gamma$ of genus $g$ with $e$
one-valent vertices. It therefore suffices to show, for a given
$\Gamma$, finiteness of the set
\[
\big\{ u_{(V,E)}\in N\,\big|\, (V,E)\in F\Gamma,\ (\Gamma,u)\text{ type
of tropical curve},\, \Delta(\Gamma,u)=\Delta\big\}.
\]
of possible slopes.

Choose a basis $\lambda_1,\ldots,\lambda_n\in \Hom(N,\ZZ)$ and define
projections $\pi_i=(\lambda_i,\lambda_n): N_\QQ\to \QQ^2$,
$i=1,\ldots,n-1$. If $h:\Gamma\to N_\RR$ is a tropical curve of degree
$\Delta$ then $\pi_i\circ h$ is a plane tropical curve of degree
$\pi_i(\Delta)$, possibly after removing contracted edges from the
domain. By \cite{mikhalkin}, Corollary~3.16 and Remark~3.17 any plane
tropical curve is dual to an integral subdivision of an integral
polygon associated to its degree. This shows that for any $i$ there
are only finitely many slopes possible for the $1$-cells of
$\pi_i\circ h(\Gamma)$. After replacing $\lambda_n$ by $\lambda_n+
\sum_{i=1}^{n-1} \eps_i\lambda_i$ for general $\eps_i\in\ZZ$, we may
assume that vertical slopes do not occur. In other words,
$\lambda_n\circ h$ is non-constant on any $E\in\Gamma^{[1]}$.

Now for non-vertical affine lines $L_i\subset \QQ^2$,
$i=1,\ldots,n-1$, the hyperplanes $\pi_i^{-1}(L_i)$ intersect
transversally and hence $\bigcap_{i=1}^{n-1} \pi_i^{-1}(L_i)\subset
N_\QQ$ is a line. Therefore
\[
\bigcap_{i=1}^{n-1} \pi_i^{-1} \big( \pi_i\circ h(\Gamma)\big)
\]
is a one-dimensional cell complex containing $h(\Gamma)$. But by what
we said before there are only finitely many slopes possible for the
edges of $\pi_i(h(\Gamma))$. Therefore only finitely many slopes may
occur in $h(\Gamma)$ for any of the tropical curves considered.
\qed

\begin{remark}
The projection method employed in the proof also shows that any
tropical curve is contained in a complete intersection of $n-1$
tropical hypersurfaces that are pull-backs from tropical curves on
$\QQ^2$ of known degree. As the types of the latter tropical curves
can be effectively enumerated by integral subdivisions of an integral
polygon this gives an algorithm to list all types of tropical curves
of given degree in any dimension. This is certainly a rather crude
method, but the only one known to the authors at the time of writing.
\end{remark}

Let us now turn to the prime transversality result for tropical
curves.

\begin{definition}\label{non-general tropically}
Let $\Delta\in\map(N\setminus\{0\},\NN)$ be a degree and
$e:=|\Delta|$. An affine constraint $\mathbf{A}=
(A_1,\ldots,A_l)$ of codimension $\mathbf{d}= (d_1,\ldots,d_l)$ is
\emph{general} for $\Delta$ if $\sum_i d_i=e+n-3$ and if for any
rational $l$-marked tropical curve $(\Gamma,\mathbf{E},h)$ of degree
$\Delta$ and matching $\mathbf{A}$ the following holds:
\begin{enumerate}
\item[(i)] $\Gamma$ is trivalent.
\item[(ii)] $h(\Gamma^{[0]})\cap \bigcup_i A_i=\emptyset$.
\item[(iii)] $h$ is injective for $n>2$. For $n=2$ it is at least
injective on the subset of vertices, and all fibers are finite.
\end{enumerate}
Otherwise it is called \emph{non-general}.
\end{definition}

\begin{proposition}\label{tropical transversality}
Let $\Delta\in\map(N\setminus \{0\},\NN)$ and  let $\mathbf A$ be an
affine constraint of codimension $\mathbf d=(d_1,\ldots,d_l)\in
\NN^l$ with $\sum_i d_i=|\Delta|+n-3$. Denote by
$\A:=\textstyle{\prod_{i=1}^l} N_\QQ/ L(A_i)$ the space of affine
constraints that are parallel to $\mathbf A$. Then the subset
\[
\mathfrak{Z}:=\Big\{ \mathbf{A}'\in \A\,\Big|\, \mathbf{A}'
\text{ is non-general for } \Delta\,\Big\}
\]
of $\A$ is nowhere dense.

Moreover, for any $\mathbf{A}'\in \A \setminus \mathfrak{Z}$ and any
$l$-marked type $(\Gamma,\mathbf{E},u)$ of genus $0$ and degree
$\Delta$ there is at most one tropical curve $(\Gamma,\mathbf{E},h)$
of type $(\Gamma,\mathbf{E},u)$ matching $\mathbf{A}'$.
\end{proposition}
\proof
(Cf. \cite{mikhalkin}, Proposition~4.11.) By Lemma~\ref{finiteness of
types} we may as well fix the type $(\Gamma,\mathbf{E},u)$ of the
(rational) tropical curves $(\Gamma,\mathbf{E},h)$ in the definition
of $\mathfrak{Z}$. Then $e:=|\Delta|$ is the number of unbounded edges
of $\Gamma$. Define the incidence variety
\[
\mathfrak{I}:=\big\{ \big((\Gamma,\mathbf{E},h),\mathbf{A}'\big) \in
\T_{(\Gamma,\mathbf{E},u)}\times \A \,\big|\, \forall i: h(E_i)\cap
A'_i\neq\emptyset\big\},
\]
and consider the diagram of canonical projections
\begin{eqnarray}\label{incidence diagram}
\begin{CD}
\mathfrak{I}@>\rho >> \T_{(\Gamma,\mathbf{E},u)}\\
@V\pi VV\\
\A.
\end{CD}
\end{eqnarray}
Let $\mathfrak{I}^\text{ng}\subset \mathfrak{I}$ be the subset of
pairs $((\Gamma,\mathbf{E},h),\mathbf{A}')$ with
$(\Gamma,\mathbf{E},h)$ non-general for $\mathbf{A}'$. We are
interested in $\mathfrak{Z}=\pi(\mathfrak{I}^\text{ng})\subset \A$ and
in the fibers of $\pi$ over the complement of this set.

We claim that the maps in Diagram~\ref{incidence diagram} can be
described by affine maps. As discussed after the dimension formula
(\ref{dim space of tropical curves}) there is an embedding of
$\T_{(\Gamma,\mathbf{E},u)}$ into $N_\QQ\times \QQ^{e-3}$ as open
subset. The factors are the position of one vertex and the lengths of
the bounded edges, respectively. Any point from $N_\QQ
\times\QQ^{e-3}$ still corresponds uniquely to a map $h:\Gamma\to
N_\QQ$, but some edges might be contracted or point in the opposite
directions than indicated by the type. Now the incidence condition of
the $i$-th edge with the $i$-th affine subspace is expressed by an
affine equation. This proves the claim. Note that our discussion
expresses $\mathfrak{I}$ as the intersection of
$N_\QQ\times\QQ_{>0}^{e-3}\times \A$ with an affine subspace. Thus
$\mathfrak{I}$ is naturally the interior of a convex polytope in a
$\QQ$-vector space.

To show that $\pi$ is either an isomorphism or has image contained in
a proper submanifold it remains to prove $\dim \mathfrak{I}\le
\dim\A$. For fixed $(\Gamma,\mathbf{E},h)\in
\T_{(\Gamma,\mathbf{E},u)}$ the set of $A'_i\in N_\QQ/L(A_i)$ with
$h(E_i)\cap A'_i\neq\emptyset$ is of codimension $\ge d_i$ inside
$N_\QQ/L(A_i)$. Thus the fiber of $\rho$ over $(\Gamma,\mathbf{E},h)$
has codimension at least $\sum_i d_i$ inside
$\{(\Gamma,\mathbf{E},h)\}\times\A$. By the dimension
formula~(\ref{dim space of tropical curves}) $(\Gamma,\mathbf{E},h)$
varies in a $(e+n-3-\operatorname{ov} (\Gamma))$-dimensional family,
and $e+n-3=\sum_i d_i$ by assumption. Hence the incidence variety on
the upper left of Diagram~\ref{incidence diagram} has at most the
dimension of $\A$, as claimed.

Finally, we look at the genericity conditions
Definition~\ref{non-general tropically},~(i)--(iii). (i) If $\Gamma$
has a vertex of valency larger than three then $\dim \T_{(\Gamma,
\mathbf{E},u)}< e+n-3$. (ii) For the subset of tropical curves with a
vertex $V\in\partial E_i$ mapping to $A_i$ the codimension of the
fiber of $\rho$ inside $\{(\Gamma, \mathbf{E},h)\}\times\A$ is
strictly larger than $\sum_i d_i$. (iii) The conditions that two
vertices of $\Gamma$ map to the same point, that the images of two
non-adjacent edges are contained in the same line or, for $n>2$, that
two non-adjacent edges intersect, are given by proper affine subspaces
inside $N_\QQ\times\QQ^{e-3}$. Finally, let $h:\Gamma\to N_\RR$ be a
trivalent tropical curve with $h(E_1)\cap h(E_2)\neq\{V\}$ for two
adjacent edges $E_1$, $E_2$. Then $u(E_1,V)=u(E_2,V)= -u(E_3,V)$ for
$E_3$ the third edge emanating from $V$, and either $h(E_1)\subset
h(E_2)$ or $h(E_2)\subset h(E_1)$, say the former. In this case the
fiber of $\pi$ through the given $((\Gamma,\mathbf E,h),\mathbf A)$ is
at least one-dimensional, for we can move the image of $V$ along the
line containing $h(E_i)$. Hence for such types $(\Gamma,\mathbf E,u)$
of tropical curves $\im(\pi)$ is contained in a proper submanifold of
$\mathfrak A$.

Thus any of these exceptional tropical curves describe subspaces of
$\mathfrak I$ with image in $\mathfrak A$ contained in a proper
submanifold. Hence $\pi(\mathfrak{I}^\text{ng}) \subset \mathfrak{A}$ is
contained in a finite union of proper submanifolds.

To summarize, two cases may occur. (1) The affine extension of $\pi$
is an isomorphism. Then $\pi(\mathfrak{I}^\text{ng}) \subset \A$ is
closed and nowhere dense. Moreover, for any
$\mathbf{A}'\in \A$ there is a unique point in $N_\QQ\times
\QQ^{e-3}$ fulfilling the extended incidence condition. This point
corresponds to a tropical curve $(\Gamma,\mathbf{E},h)\in
\T_{(\Gamma,\mathbf{E},u)}(\mathbf{A}')$ iff it lies in $N_\QQ \times
\QQ_{>0}^{e-3}$. (2) $\pi$ is not an isomorphism. Then $\im(\pi)
\subset \A$ is closed and nowhere dense, and for $\mathbf{A}'$ in the
complement of this set there is no tropical curve of type
$(\Gamma,\mathbf{E},u)$ matching $\mathbf{A}'$.
\qed
\medskip

General affine constraints imply a certain infinitesimal rigidity as
follows. Let $(\Gamma,\mathbf{E},u)$ be a type of $l$-marked rational
tropical curves and let $\mathbf{d}\in \NN^l$ with $\sum_i
d_i=|\Delta(\Gamma,u)| +n-3$. Orient the bounded edges of $\Gamma$
arbitrarily, thus defining two maps $\partial^\pm:
\Gamma^{[1]}\setminus \Gamma^{[1]}_\infty\to \Gamma^{[0]}$ with
$\partial E=\{\partial^-E,\partial^+E\}$. If $E\in\Gamma^{[1]}$ is
unbounded then $\partial^- E$ denotes the unique vertex adjacent to
$E$. For an affine constraint $\mathbf{A}$ of codimension
$\mathbf{d}$ consider the affine map
\begin{eqnarray}\label{transversality sequence}
\quad\quad\Phi: \map(\Gamma^{[0]},N_\QQ)&\lra&
\prod_{E\in\Gamma^{[1]} \setminus \Gamma_\infty^{[1]}} N_\QQ/\QQ
u_{(\partial^-E,E)} \times \prod_{i=1}^l N_\QQ/ \big(\QQ u_{(\partial^-
E_i,E_i)}+L(A_i)\big),\\
h &\longmapsto& \big((h(\partial^+E)-h(\partial^-E))_E,
(h(\partial^-E_i)-A_i)_i\big).\nonumber
\end{eqnarray}
Then $h\in\map(\Gamma^{[0]},N_\QQ)$ is the restriction to
$\Gamma^{[0]}$ of a tropical curve $\tilde h:\Gamma\to N_\RR$ of type
$(\Gamma,\mathbf{E},u)$ matching $\mathbf{A}$ iff the following
conditions are satisfied:
\begin{enumerate}
\item[(i)] $\Phi(h)=0$.
\item[(ii)] For any $E\in \Gamma^{[1]}\setminus \Gamma_\infty^{[1]}$ the
proportionality factor $\lambda\in\QQ$ with
$\partial^+E-\partial^-E=\lambda\cdot u_{(\partial^-E,E)}$ existing by (1) is
strictly positive.
\item[(iii)] For any $i$ the line passing through $h(E_i)$
intersects $A_i$ \emph{between} $h(\partial^+ E_i)$ and
$h(\partial^-E_i)$.
\end{enumerate}
The proposition now implies the following.

\begin{corollary}\label{transversality by affine map}
If the affine constraint $\mathbf{A}$ is general for $\Delta$ and
$\T_{(\Gamma,\mathbf{E},u)}(\mathbf{A})\neq\emptyset$ then the affine
map $\Phi$ in (\ref{transversality sequence}) is an isomorphism.
\end{corollary}
\proof
By the proposition there is at most one tropical curve of type
$(\Gamma,\mathbf{E},u)$ and matching $\mathbf{A}$. Thus $\Phi$ is
injective because $\Phi^{-1}(0)\neq\emptyset$ by assumption. The
proof is finished by comparing dimensions: $\Gamma$ being rational,
connected and trivalent implies $\sharp \Gamma^{[0]}= e-2$, with
$e=\sharp\Gamma^{[1]}_\infty$ the number of unbounded edges. Hence
the dimension of the left-hand side is $n(e-2)$. For the right-hand
side observe $\sharp(\Gamma^{[1]} \setminus \Gamma^{[1]}_\infty)=
\sharp\Gamma^{[0]}-1= e-3$ by rationality and connectedness, via an
Euler characteristic count. Taking into account $\sum_i d_i=e+n-3$
now also gives
\[
(e-3)(n-1)+\sum_i \big(n-(n-d_i-1+1)\big)=n(e-2).
\]
\vspace{-5ex}

\qed

\begin{remark}
In higher genus (and $n>2$) there are families of tropical curves of
larger than the expected dimension, and a result analogous to
Proposition~\ref{tropical transversality} does not hold. A treatment
of these cases therefore has to implement some kind of virtual
intersection theory on the moduli space of tropical curves. Such
spaces should itself be tropical varieties, that is, a cell complex of
integral affine polyhedra together with weights and compatible fan
structures at the vertices. The interior of each cell parametrizes
tropical curves of exactly one type. The balancing condition in higher
dimensions should come from the Minkowski weight condition formulated
in \cite{fultonsturmfels}.
\end{remark}

\section{Polyhedral decompositions and degenerations}

Here we consider degenerations of toric varieties given by toric
morphisms to $\AA^1$. They have the property that all fibers, except
the one over the origin, are isomorphic to a fixed toric variety. We
begin with some definitions.

In this paper a \emph{(rational) polyhedron} is the solution set in
$N_\QQ\simeq\QQ^n$ of finitely many linear inequalities $\langle
m,\,.\,\rangle\ge \text{const}$, $m\in M_\QQ$. As an intersection of
finitely many closed halfplanes it is always closed and convex, but
not necessarily bounded and it may have any dimension $\le n$. The
sets where some of the defining inequalities are equalities define
the \emph{faces}, which are itself lower dimensional polyhedra. A
\emph{vertex} is a zero-dimensional face. A polyhedron is
\emph{strongly convex} if it has at least one vertex.

\begin{definition}\label{polyhedral decomposition}
A (semi-infinite) \emph{polyhedral decomposition} of $N_{\QQ}$ is a
covering $\P=\{\Xi\}$ of $N_{\QQ}$ by a finite number of strongly
convex polyhedra satisfying the following properties:
\begin{enumerate}
\item[(i)]
If $\Xi \in \P$ and $\Xi' \subset \Xi$
is a face, then $\Xi' \in \P$.
\item[(ii)]
If $\Xi, \Xi' \in \P$, then $\Xi \cap \Xi'$
is a common face of $\Xi$ and $\Xi'$.
\end{enumerate}
\end{definition}
The unbounded elements of a polyhedral decomposition $\P$ define a
fan $\Sigma_\P$ by rescaling by $a\in\QQ_{>0}$ and letting $a$ tend
to $0$. In fact, for each $\Xi\in\P$ the limit $\lim_{a\to0}
a\Xi$ exists in the Hausdorff sense. Note that for every bounded
$\Xi$ this limit is just $0\in N_\QQ$. We call
\[
\Sigma_\P:=\big\{\lim_{a\to0} a\Xi\subset N_\QQ\,\big|\,
\Xi\in\P\big\}
\]
the \emph{asymptotic fan} of $\P$. The terminology is justified by
the following lemma.

\begin{lemma}
$\Sigma_\P$ is a complete fan.
\end{lemma}
\proof
If $\Xi\in\P$ is the solution set to $m_i \ge c_i$ for
$m_i\in M_\QQ$, $c_i\in\QQ$, then
\begin{eqnarray}\label{lim sigma}
\lim_{a\to 0} a\,\Xi= \big\{n\in N_\QQ\,\big|\,
\forall i\,\langle m_i,n\rangle\ge 0\big\}.
\end{eqnarray}
This is a strongly convex polyhedral cone. The intersection
of any two cones is again of this form because of (ii) in
Definition~\ref{polyhedral decomposition}. Completeness is clear.
\qed
\medskip

A polyhedral decomposition defines a degenerating family of toric
varieties as follows. First, we extend the fan $\Sigma_\P$ in $N_\QQ$
to a fan $\widetilde\Sigma_\P$ in $N_{\QQ} \times \QQ$ covering the
half-space $N_{\QQ} \times \QQ_{\ge0}$: For each $\Xi\in\P$ let
$C(\Xi)$ be the closure of the cone spanned by $\Xi \times
\{1\}$ in $N_\QQ\times\QQ$:
\[
C(\Xi) = \overline{ \big\{a \cdot (n, 1)\,\big|\, a \ge 0, n\in
\Xi\big\} }.
\]
If $\Xi$ is defined by inequalities $m_i\ge c_i$ then
\begin{eqnarray}\label{C(sigma)}
C(\Xi) = \big\{ (n,b)\in N_\QQ\times\QQ_{\ge 0}\,\big|\,
\langle m_i,n \rangle -b\cdot c_i\ge 0\big\}.
\end{eqnarray}
Thus any $C(\Xi)$ is a strongly convex polyhedral cone and
\[
\widetilde\Sigma_\P:=\big\{ \sigma\subset C(\Xi) \text{
face}\,\big|\, \Xi\in \P\big\}
\]
is a fan covering $N_\QQ\times\QQ_{\ge0}$. For later reference we
note:

\begin{lemma}\label{extended fan versus asymptotic fan}
If we identify $N_\QQ$ with $N_\QQ\times\{0\}\subset N_\QQ\times\QQ$
then
\[
\Sigma_\P= \big\{ \sigma \cap (N_\QQ\times\{0\})\,\big|\,
\sigma\in\widetilde\Sigma_\P\big\}.
\]
\end{lemma}
\proof
For $\Xi\in\P$ the description of the elements of $\Sigma_\P$ in
(\ref{lim sigma}) and of $\widetilde\Sigma_\P$ in (\ref{C(sigma)}) by
inequalities readily shows $\lim_{a\to 0} a\Xi= C(\Xi)\cap
(N_\QQ\times\{0\})$.
\qed
\medskip

The projection $N_\QQ\times\QQ\to (N_\QQ\times\QQ)/N_\QQ=\QQ$ onto
the second factor defines a non-constant map of fans from
$\widetilde\Sigma_\P$ to the fan $\{0,\QQ_{\ge0}\}$ of $\AA^1$. Thus
there is a non-constant, hence flat, toric morphism
\[
\pi: X(\widetilde\Sigma_\P)\lra \AA^1.
\]
It is equivariant for the morphism of algebraic tori
$\GG_m^{n+1}\simeq \GG(N\times\ZZ)\to \GG(\ZZ)\simeq\GG_m$. Because
$\GG(\ZZ)$ acts transitively on the closed points of $\AA^1\setminus
\{0\}$ the closed fibers of $\pi$ over $\AA^1\setminus\{0\}$ are all
pairwise isomorphic.

\begin{lemma}\label{general fiber}
For a closed point $t\in\AA^1\setminus\{0\}$ the fiber
$\pi^{-1}(t)\subset X(\widetilde\Sigma_\P)$ with the action of
$\GG(N)\subset \GG(N\times\ZZ)$ is torically isomorphic to
$X(\Sigma_\P)$.
\end{lemma}

\proof
Because $N\subset N\times\ZZ$ is the kernel of the projection to
the second factor, $\GG(N)\subset\GG(N\times\ZZ)$ induces a trivial
action on $\AA^1$, and hence respects the fibers of $\pi$.

The elements of $\widetilde\Sigma_\P$ contained in $N_\QQ\subset
N_\QQ\times\QQ$ form the fan of the open part $X(\tilde
\Sigma_\P)\setminus \pi^{-1}(0)= \pi^{-1}(\AA^1\setminus\{0\})$.
Hence Lemma~\ref{extended fan versus asymptotic fan} shows
\[
\pi^{-1}(\AA^1\setminus\{0\}) =\big(\AA^1\setminus\{0\}\big)
\times X(\Sigma_\P),
\]
with $\pi$ the projection to the first factor.
\qed
\medskip

Note also that $\widetilde\Sigma_\P$ and the product fan
$\Sigma_\P\times\{0,\QQ_{\ge0}\}$ have common refinements that agree
on $N\times\{0\}$. Subdivisions of fans lead to toric birational
morphisms. Hence there is a birational transformation
\[
X(\widetilde\Sigma_\P)\,\cdots\!\to X(\Sigma_\P)\times\AA^1,
\]
with centers on the central fiber. The central fiber can also be
easily read off from $\P$ provided that $\P$ is \emph{integral}, that
is, $\P^{[0]}\subset N$. (Integrality is a necessary and
sufficient condition for $\pi^{-1}(0)$ to be reduced.) For every
vertex $v\in\P$ the star of $\P$ at $v$ defines a complete fan
\[
\Sigma_v:=\big\{\QQ_{\ge 0}\cdot
(\Xi-v) \,\big|\, \Xi\in\P,\Xi\ni v\big\}.
\]
More generally, for $\Xi\in\P$ the rays emanating from $\Xi$ through
adjacent $\Xi'\in\P$ define a complete fan $\Sigma_\Xi$ in
$N_\QQ/L(\Xi)$:
\[
\Sigma_\Xi:=\big\{ \QQ_{\ge 0}\cdot (\Xi'-\Xi)\subset N_\QQ/L(\Xi)
\,\big|\, \Xi'\in\P, \Xi\subset\Xi' \big\}.
\]
For shortness write $X_\Xi:= X(\Sigma_\Xi)$. In particular,
$X_v=X(\Sigma_v)$ is a toric divisor in $X(\widetilde\Sigma_\P)$. For
any vertex $v\in\Xi$ let $\Sigma_v(\Xi)\subset\Sigma_v$ be the subfan
of cones intersecting $\Xi-v$ non-trivially. Then $N_\QQ\to
N_\QQ/L(\Xi)$ defines a map of fans $\Sigma_v(\Xi)\to \Sigma_\Xi$,
hence a toric morphism between the open subset
$X(\Sigma_v(\Xi))\subset X_v$ to $X_\Xi$. This map induces an
isomorphism between the $(\codim\Xi)$-dimensional toric stratum in
$X_v$ corresponding to $\Xi$ and $X_\Xi$. Thus for any $\Xi,\Xi'\in
\P$ with $\Xi'\subset\Xi$ there is a closed embedding $X_\Xi\to
X_{\Xi'}$ that is equivariant with respect to $\GG(N)$ and compatible
with compositions. While this is standard in toric geometry we wish
to view the $X_\Xi$ together with these closed embeddings as a
directed system of schemes with $\GG(N)$-action.

\begin{proposition}\label{central fiber}
Assume $\P^{[0]}\subset N$. Then there exists a system of closed
embeddings $X_\Xi\to \pi^{-1}(0)$, $\Xi\in\P$, compatible with the
directed system, inducing an isomorphism $\displaystyle
\pi^{-1}(0)\simeq \lim_{\displaystyle
\genfrac{}{}{0pt}{1}{\lra}{\Xi\in\P} } X_\Xi$.
\end{proposition}
\proof
This is a trivial special case of the discussion in \cite{logmirror},
\S2.2. (The boundedness of the cells of $\P$ that we assumed in this
paper is irrelevant for this part.) Let us just show here how to get
the standard toric embeddings $X_v\to \pi^{-1}(0)$, leaving the
straightforward verifications of compatibility with the directed
system and the universal property to the reader.

The total space $X(\widetilde\Sigma_\P)$ of the degeneration $\pi$ is
glued from affine sets $\Spec \kk[C(\Xi)^\vee \cap(M\times\ZZ)]$ for
$\Xi\in\P$. The monomial $\chi^{(0,1)}$ is the image of the affine
coordinate on $\AA^1$, so this generates the ideal of the central
fiber. On the other hand, if $\P$ is integral this equals the ideal
generated by $\chi^{(m,k)}$ for all $(m,k)\in \Int(C(\Xi)^\vee)\cap
(M\times\ZZ)$. In fact, $(m,k)\in M\times\ZZ$ lies in
$\Int(C(\Xi)^\vee)$ iff $\langle m, v\rangle +k>0$ for all $v\in
\Xi$. If all vertices of $\Xi$ are integral it suffices to test this
inequality for $v\in \Xi\cap N$, and then it holds $\langle m,
v\rangle +k\ge 1$. This implies $(m,k-1)\in C(\Xi)^\vee\cap
(M\times\ZZ)$ and hence $\chi^{(m,k)}\in (\chi^{(0,1)})$.  Thus
$\pi^{-1}(0)$ is covered by $\Spec \kk[\partial C(\Xi)^\vee
\cap(M\times\ZZ)]$, where this notation means
\[
\chi^{(m_1,a_1)}\cdot\chi^{(m_2,a_2)}:= \left\{\begin{array}{ll}
\chi^{(m_1+m_2,a_1+a_2)},&\text{if }(m_1+m_2,a_1+a_2)\in \partial
C(\Xi)^\vee,\\[1ex]
0,&\text{otherwise.}
\end{array}\right.
\]
This space has irreducible components $\Spec \kk[C_v\cap
(M\times\ZZ)]$ with $C_v\subset C(\Xi)^\vee$ the $n$-dimensional face
dual to $\QQ_{\ge 0}\cdot (v,1)$, for $v$ a vertex of $\Xi$. Now the
projection $M\times\ZZ\to M$ induces an isomorphism of $C_v\cap
(M\times\ZZ)$ with the integral points of $\big(\QQ_{\ge
0}(\Xi-v)\big)^\vee$. Hence there is a canonical isomorphism of a
toric affine patch of an irreducible component of $\pi^{-1}(0)$ with
an affine patch of $X_v$. The identification is compatible with the
gluing maps and thus gives the claimed closed embedding $X_v\to
\pi^{-1}(0)$.
\qed
\medskip

Now we describe how an affine constraint in $N_{\QQ}$ yields a family
of incidence conditions in the degeneration $X(\widetilde\Sigma_\P)
\to \AA^1$. An affine subspace $A\subset N_\QQ$ spans the linear
subspace $LC(A)\subset N_\QQ\times\QQ$ where the fan
$\widetilde\Sigma_\P$ lives. For any closed point $P$ in the big
torus of $X(\widetilde\Sigma_\P)$ the closure of the orbit
$\GG(LC(A)\cap (N\times\ZZ)).P$ defines a subvariety $Z\subset
X(\widetilde\Sigma_\P)$ projecting onto $\AA^1$. Our incidence
condition is non-trivial intersection with $Z$. Most of the relevant
properties of this subvariety generalize to any pair (fan in
$\QQ$-vector space, linear subspace). For simplicity we therefore
revert to the notation $\Sigma$ for the fan and $N_\QQ$ for the
$\QQ$-vector space during the following discussion.

\begin{proposition}\label{orbit closures for affine constraints}
Let $\Sigma$ be a fan in $N_\QQ$, $P\in X(\Sigma)$ a closed point in
the big torus and $L\subset N_\QQ$ a linear subspace. For
$\sigma\in\Sigma$ denote by $X_\sigma\subset X(\Sigma)$ the
corresponding toric subvariety. Then the following holds.
\begin{enumerate}
\item $\overline{\GG(L\cap N).P}\cap \Int(X_\sigma)=\emptyset\quad
\Longleftrightarrow \quad L\cap \Int\sigma
=\emptyset$.
\item If $L\cap \Int\sigma\neq\emptyset$ then
$\GG(L\cap N)$ acts transitively on $\overline{\GG(L\cap N).P} \cap
\Int(X_\sigma)$.
\end{enumerate}
\end{proposition}
\proof
As the statements are local along $\Int(X_\sigma)$ we may assume
$\Sigma$ is the fan of faces of $\sigma$. Then $X(\Sigma)=\Spec
\kk[\sigma^\vee\cap M]$ and $X_\sigma\simeq \GG(N/L(\sigma)\cap N)$.
In a first step we reduce to the case $L(\sigma)=N_\QQ$.

Choose a complement $Q_\QQ\subset N_\QQ$ of $L(\sigma)$ with
$L= (L\cap Q_\QQ)+(L\cap L(\sigma))$ and write $Q=Q_\QQ\cap N$.
Then the natural map
\[
Q\times\big( L(\sigma)\cap N\big)\lra N
\]
is the inclusion of a sublattice of finite index. Let $\Sigma'$ be the
fan of faces of the preimage of $\sigma$ in $Q_\QQ\times L(\sigma)$.
Then the toric morphism $X(\Sigma')\to X(\Sigma)$ is a finite
surjection that is equivariant for $\GG(Q)\times \GG(L(\sigma)\cap N) \to
\GG(N)$. Under this homomorphism the action of $\GG(L\cap Q)\times
\GG(L\cap L(\sigma)\cap N)\subset \GG(Q)\times \GG(L(\sigma)\cap N)$
covers the action of $\GG(L\cap N)\subset \GG(N)$. Hence the closure
of a $\GG(L\cap N)$-orbit is the image of the closure of a $\GG(L\cap
Q)\times \GG(L\cap L(\sigma)\cap N)$-orbit and all our statements can
be checked on $X(\Sigma')$. We may thus assume $N=Q\times\big(
L(\sigma)\cap N\big)$.

Let $\overline\Sigma$ be the fan in $L(\sigma)$ of faces of $\sigma$.
Then the splitting $N=Q\times\big( L(\sigma)\cap N\big)$ defines a
decomposition $X(\Sigma)= \GG(Q)\times X(\overline\Sigma)$ that is
compatible with the product $\GG(N)= \GG(Q)\times \GG(L(\sigma)\cap
N)$, and $X_\sigma=\GG(Q)\times 0$ for $0\subset X(\overline\Sigma)$
the unique zero-dimensional torus orbit. In particular,
$\overline{\GG(L\cap N). P} \cap X_\sigma =\emptyset$ if and only if
$\overline{\GG(L\cap L(\sigma)\cap N). \Psi(P)}\cap 0=\emptyset$ where
$\Psi: X(\Sigma)\to X(\overline\Sigma)$ is the projection onto the
second factor. Thus we may verify the conclusions of the
proposition after dividing out $\GG(Q)$. This amounts to going over
from $\Sigma$ to $\overline\Sigma$, from $N_\QQ$ to $L(\sigma)$ and
from $L$ to $L\cap L(\sigma)$. Hence we may assume $L(\sigma)=
N_\QQ$. Then (2) follows trivially from (1) because $X_\sigma=0$ is
just a closed point.

For (1) assume that $0\not\in \overline{\GG(L\cap N).P}$. Then there
exists $f=\sum_{m\in\sigma^\vee\cap M} a_m \chi^m \in\O(X(\Sigma))$
with
\[
f(0)=0,\quad f|_{\overline{\GG(L\cap N).P}}=1.
\]
The $\GG(L\cap N)$-invariant part $\sum_{m\in\sigma^\vee\cap
L^\perp\cap M} a_m \chi^m$ of $f$ has the same properties, and is
thus non-constant. Hence there exists $m\in M$ with $m|_L=0$,
$m|_{\Int(\sigma)} >0$. This implies $L\cap
{\Int(\sigma)}=\emptyset$. Conversely, if $L\cap
{\Int(\sigma)}=\emptyset$ there exists $m\in L^\perp\cap M$ with
$m|_{\Int(\sigma)} >0$, and then $f=\chi^m$ provides a $\GG(L\cap
N)$-invariant function separating $0$ from $\overline{\GG(L\cap
N).P}$.
\qed
\medskip

Note that the intersection of $\overline{\GG(L\cap N).P}$ with the
toric boundary of $X(\Sigma)$ generally need not be reduced. Consider
for example $N_\QQ=\QQ^2$, $L=\QQ\cdot(2,3)$ and $\Sigma$ the fan of
faces of $\sigma=\QQ_{\ge 0}^2$. Then $\overline{\GG(L\cap N).(1,1)}$
is the Neil parabola and the intersection with the toric boundary is
a point of multiplicity $5$. However, in the proof of the proposition
we have indeed already shown:

\begin{corollary}\label{transverse orbit closures}
In the situation of the proposition assume that $\sigma\in \Sigma$ is
contained in $L$. Let $Q$ be a complement to $L(\sigma)\cap N$ in $N$
with $L\subset (L\cap Q_\QQ)+(L\cap L(\sigma))$, and let
$\overline\Sigma$ be the fan in $L(\sigma)$ of faces of $\sigma$.
Then locally around the big torus $\GG_\sigma\subset X_\sigma$ there
is a toric isomorphism of $X(\Sigma)$ with $\GG(Q)\times
X(\overline\Sigma)$ mapping $\overline{\GG(L\cap N).P}$ to $
\GG(L\cap Q) \times \overline{\GG (L(\sigma)\cap N).P'}$ for some
$P'\in \GG_\sigma$. 
\qed
\end{corollary}

Reverting to our situation of a toric degeneration associated to a
polyhedral decomposition $\P$ we obtain the following.

\begin{corollary}\label{transversality of incidence conditions}
Let $\P$ be an integral polyhedral decomposition, $A\subset N_\QQ$ an
affine subspace and $v\in A\cap N$ a vertex of $\P$. Consider the
associated toric degeneration $\pi: X(\tilde\Sigma_\P)\to \AA^1$ with
reduced central fiber $\pi^{-1}(0)= \bigcup_{v'\in\P^{[0]}} X_{v'}$
and let $P\in \Int(X(\tilde\Sigma_P))$ be a closed point. Then locally
around $\Int(X_v)$ in $X(\tilde\Sigma_\P)$ the orbit closure
$Z=\overline{\GG(LC(A)\cap(N\times\ZZ)).P}\subset X(\tilde\Sigma_\P)$
projects smoothly to $\AA^1$ with fiber over $0\in\AA^1$ a
$\GG(L(A)\cap N)$-orbit. Moreover, the map $L(A) \otimes_\QQ \O_Z\to
\Theta_{Z/\AA^1}$ induced by the fiberwise action of $\GG(L(A)\cap
N)\subset \GG(N\times\{0\})$ on $Z\to \AA^1$ is an isomorphism.
\end{corollary}

\proof
Because $v$ is integral we have the splitting $N\times\ZZ=
(N\times\{0\})\oplus \big(\ZZ\cdot(v,1)\big)$. This induces an
isomorphism of the toric affine open set $U\simeq
\AA^1\times\Int(X_v)\subset X(\tilde\Sigma_\P)$ corresponding to
$\QQ_{\ge 0} (v,1)\in \tilde\Sigma_\P$ with $\GG(N)\times \AA^1$, and
such that $\pi$ is the projection to the second factor. By
Corollary~\ref{transverse orbit closures} the orbit closure
$\overline{\GG(LC(A) \cap(N\times\ZZ)).P}$ maps to $\GG(L(A)\cap
N).P'\times \AA^1$ under this isomorphism, for some $P'\in \Int(X_v)$.
From this the statements are evident.
\qed
\bigskip

The final topic of this section concerns refinements of polyhedral
decompositions.

\begin{proposition}\label{adapted polyhedral decomposition}
Let $h:\Gamma\to N_\RR$ be a tropical curve and $S=\{\QQ_{\ge
0}(h(E)-h(\partial E))\, |\, E\in \Gamma_{\infty}^{[1]} \}$ the set
of directions of unbounded edges. Then for any fan $\Sigma$ on
$N_\QQ$ with $S\subset\Sigma^{[1]}$ there exists a polyhedral
decomposition $\P$ of $N_\QQ$ with asymptotic fan $\Sigma$ and such
that
\[
\bigcup_{b\in\Gamma^{[\mu]}} h(b)\subset
\bigcup_{\Xi\in\P^{[\mu]}} \Xi,\quad \mu=0,1.
\]
Moreover, $\P^{[0]}$ can be chosen to contain any finite subset of
$N_\QQ$.
\end{proposition}

\proof
If $\P$ is a polyhedral decomposition of $N_\QQ$ and $\Xi\subset
N_\QQ$ is a polyhedron then $\big\{\Xi\cap\Xi' \,\big|\,
\Xi'\in\P\big\}$ is a face-fitting decomposition of $\Xi$ into
subpolyhedra. Thus if $\P$, $\P'$ are two polyhedral decompositions
of $N_\QQ$ then $\big\{\Xi\cap\Xi' \,\big|\, \Xi\in\P,
\,\Xi'\in\P'\big\}$ is a refinement of both $\P$ and $\P'$. Moreover,
if $\P$, $\P'$ have the same asymptotic fan $\Sigma$ then so does
this common refinement. It remains to construct, for every edge
$E\in\Gamma^{[1]}$ a polyhedral decomposition $\P$ with $h(E)\in
\P^{[1]}$. If $E$ is an unbounded edge with $\partial E=V$ then by
assumption $h(V)+\Sigma$ is such a decomposition. We may therefore
assume $E$ to be bounded. Denote by $V_1,V_2$ the adjacent vertices.

Let $\Sigma'\subset\Sigma$ be the subset of $n$-dimensional cones
containing the ray $e:=\QQ_{\ge 0}(h(E)-h(V_1))$ and
$B_1=\bigcup_{\sigma\in\Sigma'} \sigma$. We claim that $\partial B_1$
divides $N_\QQ$ into two connected components. In fact, $\Sigma$
consists of cones over cells of the polyhedral decomposition
$\{\sigma\cap S^{n-1}\,|\, \sigma\in \Sigma\}$ of the unit sphere in
$N_\QQ$. From this point of view $\Sigma'$ is obtained by taking
cones over the union of $(n-1)$-cells containing the one-point set
$e\cap S^{n-1}$. This is nothing but the closed star of the minimal
cell containing $e\cap S^{n-1}$, which is thus a cell. Therefore its
boundary divides $S^{n-1}$ into two connected components. Taking
cones gives the desired statement for $N_\QQ$.

We have thus divided $N_\QQ$ into the two regions $B_1$ and
$B_2=N_\QQ\setminus \Int B_1 = \bigcup_{\sigma\in\Sigma^{[n]}
\setminus\Sigma'} \sigma$ intersecting along a union of $(n-1)$-faces
of $\Sigma$. Any maximal cone of $\Sigma$ is contained in one of the
two regions, and $e\setminus\{0\}\subset\Int(B_1)$,
$-e\setminus\{0\}\subset\Int(B_2)$. Now define $\P$ as the set of
polyhedra of the following types: (I)~$h(V_1)+\sigma$,
$\sigma\in\Sigma$, $\sigma\subset B_2$ (II)~$h(V_2)+\sigma$,
$\sigma\in\Sigma$, $\sigma\subset B_1$ (III)~$h(E)+\tau$,
$\tau\in\Sigma$, $\tau\subset B_1\cap B_2$. Polyhedra of types~I and II
give polyhedral decompositions of the disjoint regions $h(V_2)+B_1$ and
$h(V_1)+B_2$, respectively. The closure of the complement of their union
is covered by polyhedra of type~III. Note that such a polyhedron
$h(E)+\tau$ has boundary
\[
\partial (h(E)+\tau)= (h(E)+\partial\tau) \cup (h(V_1)+\tau)
\cup (h(V_2)+\tau),
\]
which is thus covered by polyhedra of types~III, I, and II, in the
order of appearance. This shows that $\P$ is a polyhedral
decomposition. Moreover, in the above notation the limits under
rescaling $\lim_{a\to 0}a\cdot \Xi$ for $\Xi\in \P$ of types I, II
and III are $\sigma$, $\sigma$ and $\tau$, respectively. Hence the
asymptotic fan of $\P$ is $\Sigma$.
\qed

\section{Dual intersection graphs and tropical curves}

In the last section we saw how a polyhedral decomposition $\P$ with
asymptotic fan $\Sigma_\P$ gives rise to a toric degeneration
$X\to\AA^1$ with general fibers $X(\Sigma_\P)$. The central fiber is
a union of toric varieties, one for each vertex of $\P$, glued along
toric divisors: $X_0=\bigcup_{v\in\P^{[0]}} X_v$. The topic of this
section is the basic correspondence between tropical curves factoring
over the $1$-skeleton of $\P$ and stable maps to $X_0$ of a certain
primitive kind on each irreducible component and with image disjoint
from toric strata of codimension at least $2$. This last property is
very essential in what follows and therefore deserves a name.

\begin{definition}\label{torically transverse}
Let $X$ be a toric variety. An algebraic curve $C\subset X$ is
\emph{torically transverse} if it is disjoint from all toric strata
of codimension $>1$.

A stable map $\varphi: C\to X$ defined over a scheme $S$ is torically
transverse if the following holds for the restriction $\varphi_s$ of
$\varphi$ to every geometric point $s\to S$: $\varphi_s^{-1}(\Int X)
\subset C_s$ is dense and $\varphi_s(C_s)\subset X$ is a torically
transverse curve.
\end{definition}

Note that if $C\subset X$ is torically transverse then no irreducible
component of $C$ lies in the toric boundary; similarly, a torically
transverse stable map defined over $\kk$ does not contract any
irreducible component to the toric boundary

The tropical curve associated to an algebraic curve $C\subset X_0$
will be its dual intersection graph, viewed as subcomplex of $\P$.
For the balancing condition recall that if $X$ is a non-singular
toric variety and $u_1,\ldots,u_k\in N$ are the primitive generators
for the rays of the defining fan $\Sigma$ then
\[
H_2(X,\ZZ)=\Big\{(w_1,\ldots,w_k)\in\ZZ^k\,\Big|\,
\sum_i w_i u_i=0 \Big\}.
\]
The $w_i$ are the intersection numbers of a singular cycle with the
toric divisors $D_i$ corresponding to $u_i$. Thus a homology class
$(w_1,\ldots,w_k)\in H_2(X,\ZZ)$ with all $w_i\ge 0$ can be
viewed as the tropical curve $h:\Gamma \to N_\RR$ with edges the rays
$\QQ_{\ge 0} u_i$ of $\Sigma$, weighted by $w_i$, for all $i$ with
$w_i>0$. While for general $\Sigma$ a homological interpretation of
the right-hand side does not seem possible, one can still define
intersection numbers $w_i$ for transverse curves because a toric
variety is non-singular in codimension one. Thus the following result
suffices for our purposes.

\begin{lemma}\label{geometric balancing condition}
Let $X$ be a proper toric variety and $u_1,\ldots, u_k\in N$ the
primitive generators of the rays of the defining fan. Assume
$\varphi: C\to X$ is a torically transverse stable map from a complete,
normal algebraic curve. Then $\sum_{i=1}^k w_i u_i=0$ for
$w_i=\deg \varphi^*(D_i)$ the intersection number of $C$ with the
toric divisor corresponding to $u_i$.
\end{lemma}
\proof
It suffices to prove the claimed identity after pairing with any element
$m\in M=\Hom(N,\ZZ)$. Note that $\langle m,u_i\rangle$ is the order
of zero (or minus the pole order) of the monomial rational function
$\chi^m$ along $D_i$. Pulling back by $\varphi$ brings in a factor
$w_i$. Now $\chi^m$ being monomial all zeroes and poles are along
toric divisors. Thus the identity $\sum_i \langle m,w_i
u_i\rangle=0$ is nothing but the vanishing of the sum of orders of
zeros and poles of $\varphi^*(\chi^m)$ on the complete, non-singular
curve $C$.
\qed
\medskip

Having clarified this point we are in position to produce tropical
curves from certain stable maps to $X_0$. We need the following
refinement of Definition~\ref{torically transverse}.

\begin{definition}\label{pre-log curves}
Let $X_0=\bigcup_{v\in\P} X_v$ be the central fiber of the toric
degeneration $X\to \AA^1$ defined by an integral polyhedral
decomposition $\P$ of $N_\QQ$. A \emph{pre-log curve} on $X_0$ is a
stable map $\varphi: C\to X_0$ with the following properties:
\begin{enumerate}
\item[(i)] For any $v$ the projection $C\times_{X_0}X_v \to
X_v$ is a torically transverse stable map.
\item[(ii)] Let $P\in C$ map to the singular locus of $X_0$. Then $C$
has a node at $P$, and $\varphi$ maps the two branches $(C',P)$,
$(C'',P)$ of $C$ at $P$\ to different irreducible components
$X_{v'}$, $X_{v''}\subset X_0$. Moreover, if $w'$ is the intersection
index with the toric boundary $D'\subset X_{v'}$ of the restriction
$(C',P)\to (X_{v'},D')$, and $w''$ accordingly for $(C'',P)\to
(X_{v''},D'')$, then $w'=w''$.
\end{enumerate}
\end{definition}

Note that (ii) implies that if an irreducible component of $C$ hits
the intersection of two irreducible components $X_{v'}$, $X_{v''}
\subset X_0$ then $C$ has a node at this intersection with the two
branches entering both $X_{v'}$ and $X_{v''}$. There are no ``loose
ends'', so to speak. In the proof of the Main Theorem in
Section~\ref{section main theorem} we will see that (ii) is indeed a
necessary condition for a torically transverse curve to deform in the
family $\pi$. The necessity of such a condition in a similar context
was first pointed out by Tian \cite{tian}, and it occurs at prominent
place in Jun Li's work on Gromov-Witten theory for semistable
degenerations \cite{junli1}.

Now let us see how a pre-log curve $\varphi: C\to X_0$ gives rise
to a tropical curve.

\begin{construction}\label{open dual intersection graph}
We say two irreducible components of $C$ are
\emph{indistinguishable} if they intersect in a node \emph{not}
mapping to the singular locus of $X$.  Now define a weighted open
graph $\tilde\Gamma$ together with a map $h:\tilde \Gamma\to N_\RR$ as
follows.  Its set of vertices equals the quotient of the set of
irreducible components of $C$ modulo identification of
indistinguishable ones. If $C_V\subset C$ denotes the irreducible
component indexed by a vertex $V$ then $h(V)=v$ for the unique
$v\in\P^{[0]}$ with $\varphi(C_V)\subset X_v$ (Definition~\ref{pre-log
curves},~(i)). The set of bounded edges of $\tilde\Gamma$ is the set of
nodes of $C$, with $P_E\in C$ denoting the nodal point corresponding
to $E\in \tilde\Gamma^{[1]} \setminus \tilde\Gamma_\infty^{[1]}$. The
map $h$ identifies $E$ with the line segment joining $h(V')$, $h(V'')$
if $P_E\in C_{V'}\cap C_{V''}$. If $D\subset X_0$ denotes the union of
toric prime divisors of the $X_v$ \emph{not} contained in the
non-normal locus of $X_0$ then the set of unbounded edges is
$\varphi^{-1}(D)$. An unbounded edge $E$ labeled by $P_E\in
\varphi^{-1}(D)$ attaches to $V\in \tilde\Gamma^{[0]}$ if $P_E\in
C_V$. If $D_e\subset X_{h(V)}$, $e\in \Sigma^{[1]}$, is the toric prime
divisor with $\varphi(P_E)\subset D_e$ then $h$ maps $E$
homeomorphically to $h(V)+e\subset N_\QQ$. Finally define the weights
of the edges at a vertex $V$ by the intersection numbers of
$\varphi|_{C_V}$ with the toric prime divisors of $X_{h(V)}$. This is
well-defined by Definition~\ref{pre-log curves},~(ii).

While $\tilde\Gamma$ may have divalent vertices
Definition~\ref{pre-log curves},~(ii) assures that the two weights at
such a vertex agree. We may thus remove any divalent vertex by
joining the adjacent edges into one edge. The resulting weighted
open graph $\Gamma$ has the same topological realization as
$\tilde\Gamma$ and hence $h$ can be interpreted as a map $h:\Gamma\to
N_\QQ$. This is a tropical curve for the balancing condition holds by
Lemma~\ref{geometric balancing condition}.

Note that if $g(C)=0$ then $\Gamma$ must be a tree, and more generally
$b_1(\Gamma)\le g(C)$. In higher genus $\Gamma$ may have multiple
edges.
\end{construction}

\section{Maximally degenerate algebraic curves}\label{sect5}

The main result of this section is a description of the space of
pre-log curves with fixed associated tropical curve $h:\Gamma\to
N_\QQ$. This is a partial converse of Construction~\ref{open dual
intersection graph}. To make this work it is necessary to assume the
tropical curve to be trivalent and of the correct genus. In the
rational case this can be achieved by imposing a general affine
constraint (Definition~\ref{non-general tropically}). Then
there exist indeed only finitely many algebraic curves with given
tropical dual intersection graph, and their number can be readily
computed. 

On the side of algebraic curves trivalence amounts to requiring that
$\varphi:C\to X_0$ is of the following form componentwise.

\begin{definition}\label{def lines}
Let $X$ be a complete toric variety and $D\subset X$ the toric
boundary. A \emph{line} on $X$ is a non-constant, torically
transverse map $\varphi: \PP^1\to X$ such that $\sharp
\varphi^{-1}(D)\le 3$.
\end{definition}

A line intersects either $2$ or $3$ toric prime divisors. We refer to
these cases as \emph{divalent} and \emph{trivalent} respectively. For
the following discussion of lines we fix the following notation. Let
$u_i\in N$, $i=1,2$ or $i=1,2,3$ be the primitive generators of the
rays corresponding to the divisors being intersected, and $w_i$ the
intersection numbers with $\varphi$. Then $\sum_i w_i u_i=0$
(Lemma~\ref{geometric balancing condition}). Conversely, for any
$(\mathbf u, \mathbf w)=(u_i,w_i)\in N^a\times(\NN\setminus\{0\})^a$,
$a\in\{2,3\}$, with $u_i\in N$ primitive and $\sum_i w_i u_i=0$ there
exists a \emph{moduli space of lines} $\LL_{(\mathbf u, \mathbf w)}$
of \emph{type} $(\mathbf u, \mathbf w)$. It is an open subspace of an
appropriate space of stable maps. Clearly, $\GG(N)$ acts on
$\LL_{(\mathbf u, \mathbf w)}$ by composition with the action on $X$.
Let $E=(\sum_i \QQ u_i)\cap N$. Beware this need not be the same
as the lattice generated by the $u_i$.

\begin{lemma}\label{dimension reduction for lines}
Any line of type $(\mathbf u, \mathbf w)$ is contained
in the closure of a fiber of the canonical map $\GG(N)\to \GG(N/E)$.
The map
\[
\LL_{(\mathbf u, \mathbf w)}\lra \GG(N/E)
\]
thus defined is a morphism that is equivariant under $\GG(N)\to \GG(N/E)$.
\end{lemma}
\proof
Let $S$ be the one- or two-dimensional toric variety defined by the
complete fan in $E$ with rays $\QQ u_i$. Then up to a toric
birational transformation $X$ is a product $S\times Y$ with $Y$ a
complete toric variety defined by a fan in $(N/E)_\QQ$. Moreover, we
can take this transformation to be an isomorphism at the generic
points of the divisors $D_i\subset X$ corresponding to the $u_i$.
Then it is an isomorphism in a neighborhood of each line of the
considered type. This identifies $\LL_{(\mathbf u,
\mathbf w)}$ with the space of lines
$\tilde\varphi: \PP^1\to S\times Y$ intersecting only the strict
transforms of $D_i$. In particular, $\tilde\varphi$ is disjoint from
toric divisors of $S\times Y$ projecting onto $S$. Hence the
composition of $\tilde\varphi$ with the projection $S\times Y\to Y$
is constant. The other statements are then clear.
\qed
\medskip

Our discussion is thus reduced to $\dim X\le 2$. The next lemma deals
with the two-dimensional (trivalent) case. It is the toric analogue of
the discussion in \cite{mikhalkin}, Section~6.3. Let $\Sigma_{\PP^2}$
be the fan in $\QQ^2$ with rays $\QQ\cdot(1,0)$, $\QQ\cdot(0,1)$,
$\QQ\cdot (-1,-1)$.

\begin{lemma}\label{lines on surface}
Let $\Lambda$ be a free abelian group of rank $2$ and $(\mathbf u,
\mathbf w)= (u_i,w_i) \in
\Lambda^3\times(\NN\setminus\{0\})^3$ with $u_i$ primitive and
$\sum_i w_i u_i=0$. Let $S$ be the toric surface associated to the
complete fan $\Sigma_S$ with rays $\QQ u_i$ and denote by
$f_{(\mathbf u, \mathbf w)}: \PP^2\to S$ the
covering defined by the map of fans
\[
\Sigma_{\PP^2}\lra \Sigma_S,\quad
(a,b)\longmapsto a w_1\cdot u_1+ b w_2\cdot u_2.
\]
Then any line $\varphi:\PP^1\to S$ of type $(\mathbf u,
\mathbf w)$ is isomorphic to the composition of a linear
embedding $\PP^1\to \PP^2$ with $f_{(\mathbf u,
\mathbf w)}$.
\end{lemma}
\proof
Let $C$ be the normalization of an irreducible component of
$\PP^1\times_S\PP^2$ and $\tilde\varphi: C \to \PP^2$ the projection.
It suffices to show that $\tilde\varphi$ is the embedding of a line
in $\PP^2$ and the projection $C\to \PP^1$ has degree~$1$.

The degree of $f_{(\mathbf u, \mathbf w)}$ equals the index $\delta$
of the sublattice of $\Lambda$ generated by the $w_i u_i$. Now an
exercise in toric geometry shows that over the $i$-th toric divisor
$D_i\subset S$ the covering $f_{(\mathbf u, \mathbf w)}$ has
$\delta/w_i$ branches, each of which totally branched of order $w_i$
over $D_i$. This order agrees with the intersection index of $\varphi$
with $D_i$ at the unique point of intersection. Thus the composition
$C\to\PP^1\times_S\PP^2\to \PP^1$ is an unbranched cover, hence an
isomorphism, and there is a unique point of intersection of
$\tilde\varphi$ with the $i$-th toric prime divisor in $\PP^2$, with
intersection number $1$. Hence $\tilde\varphi$ is the embedding of a
curve of degree~$1$.
\qed

\begin{remark}
1)\ \ From the construction in the lemma it follows that in the
trivalent case a line in $X(\Sigma)$ is the normalization of a
rational curve with at most nodes. It touches each of the three
divisors $D_i$ in only one point in $\Int(D_i)$ of order $w_i$. In
particular, a trivalent line in $\PP^2$ can have arbitrary degree.\\
2)\ \ A divalent line is just a multiple cover of degree $w_1=w_2$ of
an orbit closure for $\GG(\ZZ u_i)$, intersecting $D_1$, $D_2$
transversally.
\end{remark}

The two previous lemmas imply the following result on the structure of
the space of lines of fixed type.

\begin{proposition}\label{space of lines is torsor}
Let $X$ be a toric variety. Then the action of the big torus $\GG(N)$
on the space of lines $\LL_{(\mathbf u, \mathbf w)}$ of fixed type is
transitive. In the trivalent case this action is simply transitive,
while in the divalent case the action factors over a simply
transitive action of $\GG(N/\ZZ u_1)= \GG(N/\ZZ u_2)$.
\end{proposition}

\proof
Lemma~\ref{dimension reduction for lines} reduces to the case $\dim
X\le 2$. In the one-dimensional (divalent) case there is only one
isomorphism class of stable maps $\PP^1\to X$ that is totally
branched over $0$ and $\infty$. Its image is the closure of
a $\GG(\ZZ u_i)$-orbit.

In the trivalent case the toric birational morphism $X\to S$ to the
toric surface from Lemma~\ref{lines on surface} reduces further to
the case $X=S$. We retain the notations from this lemma. Any line
$\varphi:\PP^1\to S$ lifts to at most $\delta=\deg f_{(\mathbf u,
\mathbf w)}$ distinct lines in $\PP^2$. On the other hand, $\delta$
is also the order of the kernel of the homomorphism $\GG(\ZZ^2)\to
\GG(\Lambda)$ underlying $f_{(\mathbf u, \mathbf w)}$, and
$\GG(\ZZ^2)$ acts transitively on the set of lines in $\PP^2$. Hence
$\ker \big(\GG(\ZZ^2)\to \GG(\Lambda)\big)$ acts transitively on the
set of lifts of $\varphi$. This proves that the action of
$\GG(\Lambda)$ on the space of lines on $S$ of considered type is
simply transitive.
\qed
\medskip

Having understood lines on toric varieties we are now in position to
discuss maximally degenerate curves on unions of toric varieties.

\begin{definition}\label{maximally degenerate curves}
Let $X_0=\bigcup_{v\in\P^{[0]}} X_v$ be the central fiber of a toric
degeneration defined by an integral polyhedral decomposition $\P$. A
pre-log curve (Definition~\ref{pre-log curves}) $\varphi: C\to X_0$
is called \emph{maximally degenerate} if for any $v\in\P^{[0]}$ the
projection $C\times_{X_0}X_v \to X_v$ is a line, or, for $n=2$, the
disjoint union of two divalent lines intersecting disjoint toric
divisors.
\end{definition}

Thus a maximally degenerate curve is a collection of lines, at most
one ($n=2$: two) for each irreducible component of $X_0$, which match
in the sense that they glue to a pre-log curve. The matching
condition involves both incidence of the intersections with the toric
divisors and equality of the intersection numbers for glued branches.

To realize a rational tropical curve $h:\Gamma\to N_\RR$ as the
dual intersection graph of a maximally degenerate curve we assume
that $\Gamma$ factors over the $1$-skeleton of a polyhedral
decomposition.

\begin{proposition}\label{count of maximally degenerate curves}
{\rm 1)}\ \ Let $\Delta\in \map(N\setminus\{0\},\NN)$ be a degree and
$\mathbf{A}=(A_1,\ldots, A_l)$ an affine constraint that is general
for $\Delta$. If $\T_{(\Gamma,\mathbf E,u)}(\mathbf{A})\neq\emptyset$
for an $l$-marked tree $(\Gamma,\mathbf E)$ then the map
\begin{eqnarray}\label{lattice map}
\map(\Gamma^{[0]},N)&\lra& \prod_{E\in\Gamma^{[1]}
\setminus \Gamma_\infty^{[1]}} N/\ZZ u_{(\partial^- E, E)} \times
\prod_{i=1}^l N/ \big(\QQ u_{(\partial^- E_i,E_i)}+L(A_i)\big)\cap N,\\
h &\longmapsto& \big((h(\partial^+E)-h(\partial^-E))_E,
(h(\partial^-E_i))_i\big).\nonumber
\end{eqnarray}
is an inclusion of lattices of finite index $\mathfrak D$. Here
$\partial^\pm: \Gamma^{[1]}\setminus \Gamma^{[1]}_\infty\to
\Gamma^{[0]}$ is an arbitrarily chosen orientation of the bounded
edges, that is, $\partial E= \{\partial^- E,\partial^+ E\}$ for any
$E\in \Gamma^{[1]}\setminus \Gamma^{[1]}_\infty$. If $E\in
\Gamma^{[1]}_\infty$ then $\partial^-E$ denotes the unique vertex
adjacent to $E$.
\smallskip

\noindent
{\rm 2)}\ \ Assume in addition that $\P$ is an integral polyhedral
decomposition of $N_\QQ$ with
\[
h(\Gamma^{[\mu]})\subset
\bigcup_{\Xi\in\P^{[\mu]}} \Xi,\quad \mu=0,1,\quad\text{ and }\quad
h(\Gamma)\cap A_j\subset \P^{[0]},\quad j=1,\ldots,l,
\]
and associated toric degeneration $X(\widetilde\Sigma_\P)\to \AA^1$
with central fiber $X_0$, and let $P_j$, $j=1,\ldots,l$ be closed
points in the big torus of $X(\widetilde\Sigma_\P)$. Then $\mathfrak
D$ equals the number of isomorphism classes of maximally degenerate
curves in $X_0$ with associated tropical curve $\Gamma$ and
intersecting
\[
Z_i:=\overline{\GG (LC(A_i)).P_i}\subset X(\widetilde\Sigma_\P),
\]
the closure of the orbit through $P_i$ for the subgroup
$\GG(LC(A_i))\subset \GG(N\times \ZZ)$ acting on $X(\widetilde\Sigma_\P)$.
\end{proposition}

\proof
Tensored with $\QQ$ the map agrees with the linear part of $\Phi$ in
Corollary~\ref{transversality by affine map}. By this corollary and
since $\mathbf{A}$ is general $\Phi$ is an isomorphism. This
proves~(1).

Going over to algebraic tori thus leads to a covering of degree
$\mathfrak D$. Now a maximally degenerate curve with associated
tropical curve $\Gamma$ is glued from a number of lines. To make this
precise let $\tilde \Gamma$ denote the open graph obtained from
$\Gamma$ by inserting vertices at all points of
$h^{-1}(\P^{[0]})\setminus \Gamma^{[0]}$. Conversely, $\Gamma$ is
obtained from $\tilde\Gamma$ by merging all edges meeting in a
divalent vertex (cf.\ Construction~\ref{open dual intersection
graph}). Now for each vertex $V\in \tilde\Gamma^{[0]}$ there is a
moduli space of lines $\LL_{(\mathbf u(V), \mathbf w(V))}$ on the
component $X_{h(V)}\subset X_0$. Here $\mathbf u(V), \mathbf w(V)$ are
the tuples of directions $u_{(V,E)}$ and weights $w(E)$ for edges $E$
adjacent to $V$, respectively. A trivalent vertex $V$ corresponds to a
trivalent line, and in this case $\LL_{(\mathbf u(V), \mathbf w(V))}$
is a torsor under $\GG(N)$ (Proposition~\ref{space of lines is
torsor}). Hence $\prod_{V\in \Gamma^{[0]}} \LL_{(\mathbf u(V), \mathbf
w(V))}$ is a torsor under $\map( \Gamma^{[0]}, \GG(N))$. If $V$ is
divalent with adjacent edges $E^\pm(V)$ then $\LL_{(\mathbf u(V),
\mathbf w(V))}$ is only a torsor under $\GG(N/\ZZ u_{(V,E^-(V))})=
\GG(N/\ZZ u_{(V,E^+(V))})$. Taken together $\prod_{V\in
\tilde\Gamma^{[0]}} \LL_{(\mathbf u(V), \mathbf w(V))}$ is a torsor
under the algebraic group belonging to the left-hand side of
$(\ref{lattice map})$ times $\GG:=\prod_{V\in \tilde\Gamma^{[0]}
\setminus \Gamma^{[0]}} \GG(N/\ZZ u_{(V,E^-(V))})$. 

Now the matching and incidence conditions are also a torsor, but
under the algebraic torus action defined by the right-hand side of
$(\ref{lattice map})$, times $\GG$. In fact, if a bounded edge $E$ of
$\Gamma$ connects two trivalent vertices then the lines $L^-\in
\LL_{(\mathbf u(\partial^- E), \mathbf w(\partial^-E))}$, $L^+\in
\LL_{(\mathbf u(\partial^+ E), \mathbf w(\partial^+E))}$ intersect
the $(n-1)$-dimensional toric variety $X_{h(E)}$ in closed points
$P^\pm$ in its big torus. This big torus is a torsor under $\GG(N/\ZZ
u_{(\partial^- E, E)}) = \GG(N/\ZZ u_{(\partial^+ E, E)})$, and
$L^\pm$ glue iff $P^-=P^+$. The gluing on the schem-theoretic level
works by proposiion~\ref{central fiber}. The other case is that one of the
vertices is divalent, say $\partial^- E$. Then $\LL_{(\mathbf
u(\partial^- E), \mathbf w(\partial^-E))}$ is a torsor under
$\GG(N/\ZZ u_{(\partial^- E,E)})\subset \GG$ and hence $L^-$ is
determined uniquely by $L^- \cap X_{h(E)}$. If $E'$ is the other edge
adjacent to $\partial^- E$ then this has the effect of transferring
the matching condition on $X_{h(E)}$ to $X_{h(E')}$. In other words,
the matchings at divalent vertices form a torsor under $\GG$.

Finally, if $h(E_i)$ intersects $A_i$ in $V\in \tilde\Gamma^{[0]}$
then the image of any $L^-\in \LL_{(\mathbf u(V), \mathbf w(V))}$ is
the closure of a $\GG(\ZZ u_{(\partial^- E_i,E_i)})$-orbit, while  by
Corollary~\ref{transversality of incidence conditions} $Z_i \cap
X_{h(V)}$ is the closure of a $\GG\big(L(A_i)\cap N)$-orbit in $\Int
(X_{h(V)})$. Choose a closed point $P_i\in Z_i\cap\Int(X_{h(V)})$,
and put $L_0^-:=\overline{\GG(\ZZ u_{\partial^- E_i,E_i)}). P_i}$.
Then for any $\kk$-rational point $\lambda$ of $\GG(N)$
\[
\lambda . L_0^-\cap Z_i\neq\emptyset
\quad\Longleftrightarrow\quad
\lambda\in\GG\big((\QQ u_{(\partial^- E_i,E_i)}+L(A_i))\cap N\big)(\kk).
\]
In fact, by the definition of $L_0^-$, $\lambda.L_0^-\cap
Z_i\neq\emptyset$ iff there exists $\lambda'\in\GG (\ZZ u_{(\partial^-
E_i, E_i)})(\kk)$ with $\lambda\lambda' P_i\in Z_i$. Since
$Z_i=\overline{\GG(L(A_i)\cap N).P_i}$ it follows that $\lambda\lambda'$ is
a $\kk$-rational point of $\GG(L(A_i)\cap N)$. Thus $\lambda$ is a
$\kk$-rational point of the subtorus of $\GG(N)$ generated by $\GG
(\ZZ u_{(\partial^- E_i,E_i)})$ and $\GG(L(A_i)\cap N)$, which is the
claimed $\GG\big((\QQ u_{(\partial^- E_i,E_i)}+L(A_i))\cap N\big)$.
Conversely, for any $\kk$-rational point $\lambda$ of this subtorus
$\lambda L_0^-\cap Z_i\neq\emptyset$. Hence $\GG\big( N/ (\QQ
u_{(\partial^- E_i,E_i)}+L(A_i))\cap N\big)$ acts simply transitively
on the incidence condition with $Z_i$.

Taken together the problem of finding maximally degenerate curves of
the considered type is governed by a map of torsors that is
equivariant for a map of algebraic tori of degree $\mathfrak D$. Hence
there are exactly $\mathfrak D$ solutions.
\qed
\medskip

\begin{remark}\label{maximal degenerate curves remark}
1)\ \ The points $P_i$ in the proposition need not be general. For
example, they could all be equal to the distinguished point
$1\in\GG(N\times \ZZ)\subset X(\widetilde\Sigma_\P)$. The point is
that we only care about genericity of the degeneration of the
$\GG(LC(A_i))$-orbits $Z_i$. This is achieved by genericity of the
incidence conditions $A_i$.\\
2)\ \ The proof also determines the number of intersection points
between $Z_i$ and the image $C_0\subset X_0$ of the constructed
maximally degenerate curves that intersect $Z_i$. In fact, they
intersect only in the big torus of one component of $X_0$, namely
$X_v$ if $A_i\cap \im(h)=\{v\}$. The proof establishes a bijection
between $Z_i\cap C_0$ and the intersection of the two subtori $\GG(\ZZ
u_{(\partial^- E_i,E_i)})$ and $\GG(L(A_i)\cap N)$ in $\GG(N)$. This
latter number of intersection points is the covering degree of
\[
\GG(\ZZ u_{(\partial^- E_i,E_i)})\times \GG(L(A_i)\cap N)
\lra \GG\big((\QQ u_{(\partial^- E_i,E_i)}+L(A_i))\cap N\big),
\]
which equals the index
\[
[\ZZ u_{(\partial^- E_i,E_i)} + L(A_i)\cap N:(\QQ u_{(\partial^-
E_i,E_i)}+L(A_i))\cap N].
\]
3)\ \ In dimension two a tropical curve may have two edges with
interiors of their images intersecting. This corresponds to two
divalent lines on the same component of $X_0$ with intersecting images.
However, as we construct stable maps $C_0\to X_0$ rather than
embedded curves intersection points of such lines are not images of
nodes of $C_0$. Hence this phenomenon is irrelevant to our treatment.
\qed
\end{remark}

\section{Transversality of degenerating families}

In the previous section we described rational curves in degenerate
toric varieties that are transverse with respect to the toric strata.
In this section we show that for any finite number of one-parameter
degenerating families of curves we can always achieve this kind of
transversality by toric blow-up of the central fiber, possibly after
base change.

\begin{lemma}\label{hypersurface transversality}
Let $X$ be a toric variety and $W\subset X$ a proper subset without
irreducible components contained in the toric boundary. Then there
exists a toric blow-up $\Upsilon:\tilde X\to X$ such that the strict
transform $\tilde W$ of $W$ under $\Upsilon$ does not contain any
$0$-dimensional toric stratum.
\end{lemma}
\proof
Let $\Sigma$ be the fan defining $X$. Any toric blow-up is defined by
a subdivision of $\Sigma$. The conclusion of the theorem is stable
under further toric blow-up, and combinatorially this translates into
a freedom of further subdivision. Now subdivisions of subfans of
$\Sigma$ are subfans of a common subdivision of $\Sigma$. We may
therefore assume that $\Sigma$ is the fan defined by the faces of a
single $n$-dimensional cone $\sigma\subset N_\QQ$. (If $\dim\sigma<n$
there are no $0$-dimensional toric strata.) Then $X=\Spec
\kk[\sigma^\vee\cap M]$ is affine and taking a general element of the
ideal defining $W$ reduces to the case that $W$ is a hypersurface.

Let $f=\sum_{p\in \sigma^\vee\cap M} a_p x^p\in \kk[\sigma^\vee\cap
M]$ be the regular function defining $W$. Denote by
\[
\Delta_f=\operatorname{conv}\Big(\bigcup_{\{p\in \sigma^\vee\cap
M|a_p\neq0\}} \{p\}+(\sigma^\vee\cap M)\Big) \subset M_\QQ
\]
the Newton polyhedron of $f$. Let $\widetilde\Sigma$ be the normal
fan of $\Delta_f$. This fan has rays $\QQ\cdot u$ with $\langle
u,p\rangle\ge 0$ for all $p\in \sigma^\vee\cap M$ and such that $u^\perp$ is parallel
to a facet of $\Delta_f$. In particular, it is a subdivision of
$\Sigma$, and hence it defines a toric blow-up $\Upsilon:\tilde X\to X$.
We claim that $\Upsilon$ has the requested properties.

Let $x\in \tilde X$ be a $0$-dimensional toric stratum and $\tau\in
\widetilde\Sigma$ be the corresponding $n$-dimensional cone. Then
$\tau$ defines the affine toric neighborhood $\Spec \kk[\tau^\vee\cap
M]$ of $x$, and $\Upsilon$ is given locally by the inclusion of rings
$\phi:\kk[\sigma^\vee\cap M]\to \kk[\tau^\vee\cap M]$. Since
$\widetilde\Sigma$ is the normal fan of $\Delta_f$ there exists a
$p_\tau\in \sigma^\vee\cap M$ with $a_{p_\tau}\neq 0$ and
$\Delta_f\subset p_\tau+ (\tau^\vee\cap M)$. In particular, every
monomial in $\phi(f)$ is divisible by $\chi^{p_\tau}$:
\[
\phi(f)=\chi^{p_\tau}\cdot
\sum_{p\in \sigma^\vee\cap M\subset\tau^\vee\cap M}
a_p \chi^{p-p_\tau}.
\]
Then $f_\tau= \sum_{p\in \sigma^\vee\cap M} a_{p_\tau +p} \chi^p$ is a
function vanishing on the strict transform $\tilde W$ of $W$. This
last argument uses the assumption that $W$ is not contained in the
toric boundary. But $f_\tau(x)= a_{p_\tau}\neq0 \in \O_{\tilde
X,x}/\maxid_x=\kk$, and hence $x\not\in \tilde W$.
\qed

\begin{proposition}\label{transversality in any dimensions}
Let $X$ be a toric variety and $W\subset X$ a closed subset of
codimension $>c$. We assume that no component of $W$ is
contained in the toric boundary. Then there exists a toric blow-up
$\Upsilon:\tilde X\to X$ such that the strict transform $\tilde W$ of $W$
under $\Upsilon$ is disjoint from any toric stratum of dimension $\le c$.
\end{proposition}

\proof
By induction on $c$ we may assume that $W$ is already disjoint from
toric strata of dimension less than $c$. Any $c$-dimensional torus
orbit in $X$ has a toric neighborhood of the form $\Spec
\kk[\tau^\vee\cap M]$ for $\tau\subset N_\QQ$ an $(n-c)$-dimensional
strongly convex cone. As in the proof of the lemma the existence of
common refinements of partial subdivisions thus reduces to the case
that the fan $\Sigma$ defining $X$ is the fan of faces of such a
$\tau$. Let $\bar\Sigma$ be the fan in $L(\tau)$ of faces of $\tau$.
A linear projection $\pi: N_\QQ\to L(\tau)$ defines a map of fans
$\Sigma\to \bar\Sigma$. Let $\Psi: X\to Y$ be the corresponding toric
morphism. A subdivision $\Sigma_\tau$ of $\bar\Sigma$ induces the
subdivision
\[
\pi_\tau^{-1}\Sigma_\tau:= \big\{\tau\cap
\pi_\tau^{-1}(\omega)\,\big|\, \omega\in\Sigma_\tau\big\}
\]
of $\Sigma$. We are going to construct $\Sigma_\tau$ such that the
strict transform of $W$ for the toric blow-up $\tilde X\to X$ defined by
this subdivision is disjoint from all $c$-dimensional toric strata
of $\tilde X$.

Now the unique $c$-dimensional toric stratum of $X$ is mapped to the
$0$-dimensional toric stratum of $Y$. Because $\dim\Psi(W) \le \dim
W<n-c=\dim Y$ and because no irreducible component of $W$ is
contained in the toric boundary of $X$, $\Psi(W)\subset Y$ is a
proper subset without irreducible components in the toric boundary of
$Y$. Thus Lemma~\ref{hypersurface transversality} provides a toric
blow-up of $Y$ such that the strict transform of $\Psi(W)$ does not
contain $0$-dimensional toric strata. The desired fan $\Sigma_\tau$
is the fan declaring this toric blow-up.
\qed

\begin{proposition}\label{stable reduction}
Let $X\to \AA^1$ be a toric degeneration with special fiber
$X_0\subset X$, let $R$ be a discrete valuation ring with residue
field $\kk$ and quotient field $K$ and $\Spec R\to \AA^1$ a dominant
morphism mapping the closed point to $0\in \AA^1$. Let $(C^*\to
X\setminus X_0, (x_1^*,\ldots,x_k^*))$ be a torically transverse
stable map with $l$ marked points $x_1^*,\ldots,x_l^*:\Spec K\to C^*$,
defined over $\Spec K$.

Then possibly after base change $\AA^1\to \AA^1$, $t\mapsto t^b$ there
exists a toric blow-up $\tilde X\to \AA^1$ with centers in $X_0$ with
the following property: $C^*$ extends to a stable map $(C\to \tilde
X,(x_1,\ldots,x_l))$ over $\Spec R$ such that for every irreducible
component $\tilde X_v\subset \tilde X_0$ the projection
$C\times_{\tilde X_0}\tilde X_v \to \tilde X_v$ is a torically
transverse stable map.
\end{proposition}
\proof
Let $\hat C^*\to C^*$ be the normalization of $C$. By Abhyankar's
Lemma (\cite{sga1}, Expos\'e~X, Lemma~3.6) there exists a totally
ramified base change $t\mapsto t^b$ such that the conductor locus of
this normalization is the image of additional sections $y_1^*,\ldots,
y_l^*$. Then $C^*$ is obtained by pairwise identification of the images
of these sections. Let $\tilde x_i^*$ be the pull-back of $x_i^*$ to $\hat
C^*$. Then for any irreducible component $\hat C^*_\mu\subset\hat C^*$
the composition $\hat C^*_\mu\to \hat C^*\to C^*\to X$ marked by those
$\tilde x_i^*$ and $y_j^*$ with image contained in $\hat C^*_\mu$ is
also a marked stable map fulfilling the hypothesis. If any of these
stable maps extend after a base change then we obtain an extension of
the original stable map by pairwise identification of the extensions
$y_j$ of $y_j^*$. This argument reduces to the case $C^*$ irreducible.

Let $W$ be the closure of the image of $C^*\to X$. Then $W$ is one- or
two-dimensional. In the one-dimensional case, by
Proposition~\ref{transversality in any dimensions}, we may assume that
$W$ is disjoint from any toric stratum of $X_0$ of dimension less than
$n=\dim X_0=\dim X-1$. Then $W\to \AA^1$ is a proper and dominant map
of curves, hence a finite surjection. The geometric fibers of
$C^*\to\Spec K$ being connected, $C^*\to W$ thus factors over $\Spec
K$ (apply Stein factorization to $C^*\to \AA^1$). It therefore just
remains to extend $(C^*\to\Spec K,(x_1^*,\ldots,x_l^*))$ as a marked
stable curve. This is possible by \cite{knudsen} after another totally
ramified base change.

In the two-dimensional case Proposition~\ref{transversality in any
dimensions} merely achieves that $W$ is disjoint from toric strata of
dimension less than $n-1$. Let $X_\tau\subset X$ be an irreducible
component of $X_0$ intersecting $W$. There is a neighbourhood
$U_\tau\subset X$ of $\Int (X_\tau)\simeq \GG_m^{n-1}$
(non-canonically) isomorphic to $\Int(X_\tau)\times V_e$
($e=e(\tau)$) with
\[
V_k =\Spec \kk[x,y,t]/(xy-t^e).
\]
Let $\tilde U^*_\tau$ be the preimage of $U_\tau$ in $C^*$.  Consider
the composition of $\tilde U^*_\tau\to U_\tau$ with the projection
$U_\tau\to V_e$. This is a dominant map of reduced $\kk$-schemes of
the same dimension, hence it is \'etale away from a nowhere dense
closed subset. Let $Z_\tau\subset V_e$ be this closed subset union the
closure of the images of the marked points, intersected with $V_e$. By
another application of Proposition~\ref{transversality in any
dimensions} we may assume that $Z_\tau$ does not contain the
zero-dimensional torus orbit of $V_e$. Do this for any $\tau$ with
$X_\tau\cap W\neq\emptyset$.

Now apply the stable reduction theorem for stable maps
(\cite{fultonpandh}, Proposition~6). This gives, after another base
change, an extension of $(C^*\to X,(x_1^*,\ldots,x_l^*))$ to a marked
stable map $(\varphi:C\to X, (x_1,\ldots,x_l))$ over $\Spec R$. By
construction the image of this map is a family of torically transverse
curves. It remains to show that the restriction $\varphi_0: C_0\to X$
to the closed fiber does not contract any component to the codimension
one locus of $X_0$. It suffices to check this statement on $\tilde
U_\tau=\varphi^{-1}(U_\tau)\subset C^*$. Thus we only have to verify
that no irreducible component of $C_0\cap \tilde U_\tau$ maps to the
singular point of $V_e$. To this end let $\tilde Z_\tau=
\varphi^{-1}(Z_\tau)$. The restriction of $\varphi$ to $\tilde
U_\tau\setminus \tilde Z_\tau$ factors over a proper map to $\Spec
(\O_{\AA^1,0})\times_{\AA^1} (U_\tau\setminus Z_\tau)\subset X$. This
map is finite except possibly over a finite set $T\subset X_0$. Hence
Stein factorization produces another extension of $\varphi|_{\tilde
U_\tau^*\setminus \tilde Z_\tau}$ over $\Spec R$. This extension glues
with $\varphi|_{C\setminus \varphi^{-1}(T)}$ to a map $\varphi': C'\to
X$ because $\varphi$ is already finite on the gluing region.
Moreover, since $Z_\tau$ contains the image of the marked points,
$\varphi'$ together with the compositions $\Spec R\stackrel{x_i}{\lra}
C\to C'$ as marked points is a stable map. By construction $C\to C'$
contracts all components of $C\cap \tilde U_\tau$ that are mapped to
the singular point of $V_e$ under the composition $\tilde U_\tau\to
U_\tau\to V_e$. Thus by the uniqueness of stable reduction we must
have $C=C'$ and there were really no such contracted components. Thus
$(\varphi:C\to X,(x_1,\ldots,x_l))$ is torically transverse.
\qed

\section{Deformation theory of maximally degenerate curves}
\label{section deformation theory}

The final ingredient in our comparison of counting of rational curves
on a toric variety versus counting of tropical curves involves
smoothing of a maximally degenerate curve $C_0\to X_0$. The crucial
technical tool for this is the deformation theory of abstract
log (-arithmic) spaces, see \cite{K.Kato}, \cite{F.Kato}. These
references also provide essentially self-contained introductions into
the subject. Section~3.1 of \cite{logmirror} also contains a quick
introduction into log geometry, covering most aspects that
concern us here. We therefore use freely the notions of this
theory. All log schemes are fine, saturated and integral. To use
results from \cite{F.Kato2} we need to work with \'etale log
structures, that is, we need to go over to the \'etale site. However,
by abuse of notation we use the same symbol for a coherent sheaf
$\shF$ in the Zariski topology and its pull-back to the \'etale
site, usually denoted $\shF^\text{\'et}$. Accordingly, the stalk of
$\shF^\text{\'et}$ over the strict henselian closure $\bar P$ of
$P\in X$ is denoted $\shF_{\bar P}$. Just for this section unadorned
letters denote log spaces or morphisms of log spaces, while
underlined letters denote the (underlying) spaces or morphisms of
schemes. As index for the structure sheaf we allow both notations.

Throughout this section let $\P$ be an integral polyhedral
decomposition of $N_\QQ$, $\underline\pi:\underline X\to
\underline\AA^1$ the associated toric degeneration and $\underline
X_0\subset \underline X$ the central fiber. Let
$(\Gamma,\mathbf{E},h)$ be an $l$-marked rational tropical curve with
image contained in the $1$-skeleton of $\P$ and matching an affine
constraint $\mathbf{A}=(A_1,\ldots,A_l)$. We assume that
$(\Gamma,\mathbf{E},h)$ and $\mathbf{A}$ are transverse in the sense
that the map $\Phi$ in (\ref{transversality sequence}) before
Corollary~\ref{transversality by affine map} is an isomorphism. We
also assume that all points of intersection of the $A_i$ with
$h(\Gamma)$ are vertices of $\P$. Let $(\underline C_0,\mathbf
x_0,\underline \varphi_0)$ be a maximally degenerate, $l$-marked,
rational stable map with dual intersection graph
$(\Gamma,\mathbf{E},h)$ and with $\underline\varphi_0(x_i)\in
\underline Z_i=\overline{\GG (LC(A_i)).P_i}\subset \underline X$ for
$j=1,\ldots,l$, as given by Proposition~\ref{count of maximally
degenerate curves}. By Corollary~\ref{transversality of incidence
conditions} the map $\underline Z_i\to \underline\AA^1$ is smooth at
$\underline \varphi_0(x_i)$ since we assumed that the unique point of
intersection of $h(E_i)$ with $A_i$ is a vertex of the polyhedral
decomposition.

On any toric variety $\underline X$ the inclusion of the toric
boundary defines a natural log structure $\shM_X\to \O_X$, and a toric
morphism is naturally a (log-) smooth morphism for these log
structures. Thus our degeneration defines a smooth morphism $\pi:X\to
\AA^1$ of log spaces. Restricting to the central fiber now defines a
smooth morphism to the standard log point
\[
\pi_0: X_0 \lra O_0:=(\Spec\kk, \NN\times\kk^\times).
\]
Recall also that the sheaf of logarithmic derivations
$\Theta_{X(\Sigma)}$ of a toric variety defined by a fan
$\Sigma$ in $N_\QQ$ is canonically isomorphic to $N\otimes_\ZZ
\O_{X(\Sigma)}$. Hence
\[
\Theta_{X/\AA^1}= (N\oplus \ZZ)\otimes_\ZZ \O_X
/(0\oplus\O_X)= N\otimes_\ZZ\O_X,
\]
and then also $\Theta_{X_0/O_0}= N\otimes_\ZZ \O_{X_0}$ by
base change.

Next we want to lift $\underline \varphi_0:\underline C_0\to
\underline X_0$ to a log morphism. Consider the diagonal map $\NN\to
\NN^2$ and the homothety $\NN\stackrel{\cdot e}{\to}\NN$ and denote by
$S_e= \NN^2\oplus_\NN\NN$ the monoid defined by push-out. Then $S_e$
has generators $a=((1,0),0)$, $b=((0,1),0)$, $c=((0,0),1)$ with single
relation $a+b=e\cdot c$, and hence $\kk[S_e]\simeq
\kk[x,y,t]/(xy-t^e)$. Recall from Section~1 the notion of total marked
weight $w(\Gamma,\mathbf{E})= \prod_{E\in\Gamma^{[1]} \setminus
\Gamma^{[1]}_\infty} w(E) \cdot\prod_{i=1}^l w(E_i)$, and from the
proof of Proposition~\ref{count of maximally degenerate curves} the
subdivision $\tilde\Gamma$ of $\Gamma$ with $\tilde\Gamma^{[0]}=
h^{-1}(\P^{[0]})$.

\begin{proposition}\label{log structures on degenerate curves}
Assume that for every bounded edge $E\subset \tilde\Gamma^{[1]}$ the
integral length of $h(E)$ is a multiple of its weight $w(E)$.
Then there are exactly $w(\Gamma, \mathbf{E})$ pairwise
non-isomorphic pairs $[\varphi_0: C_0\to X_0,\mathbf{x}_0]$
with underlying marked stable map isomorphic to
$(\underline{C_0},\mathbf{x}_0,\underline \varphi_0)$, with $\varphi_0$
strict wherever $X_0\to O_0$ is strict, and such that the
compositions $C_0\to X_0\to O_0$ are smooth and integral.
\end{proposition}
\proof
The fact that $\underline\varphi_0$ lifts to a morphism of log spaces
means that the pull-back log structure
$\underline\varphi_0^*\shM_{X_0}\to \O_{C_0}$ factors over  a log
structure $\shM_{C_0}\to \O_{C_0}$ that is smooth over $\O_0$. The
locus of strictness of the  morphism $X_0\to O_0$ is the union of the
big tori of the components of $X_0$. On this part there is no choice
because strictness says $\underline\varphi_0^*\shM_{X_0}=\shM_{C_0}$,
and this log structure is also smooth over $O_0$ trivially. The only
non-nodal points $P\in C_0$ that do not map to this locus map to
$\bigcup_{E\in\Gamma^{[1]}_\infty} X_{h(E)}$, the degeneration of the
toric boundary of the general fibers. On this part toric
transversality implies that $\underline\varphi_0^*\shM_{X_0,P}\to
\O_{C_0,P}$ has precisely one smooth extension, namely the sum of
$(\underline\pi\circ\underline\varphi_0)^* \shM_{O_0}$ and the log
structure associated to the toric divisor $X_{h(E)}\subset
X_{h(\partial E)}$, see Proposition~1.1 in \cite{F.Kato2}.

Next we show that at a node $P\in C_0$ corresponding to a bounded
edge $E\in \tilde\Gamma^{[1]}$ there are precisely $\mu:=w(E)$
pairwise non-isomorphic extensions of $\underline
\varphi_0^*\shM_{X_0}$ to a smooth structure $\shM_{C_0}$ over $O_0$.
Similar statements occur at various places in the literature, see
e.g.\ \cite{F.Kato2}, \cite{junli2}, \cite{mochizuki}, \cite{wewers}.
Due to its central importance in this paper and to avoid a technical
discussion on how it relates to known formulations, we provide full
details here.

Denote by $t$ a linear coordinate on $\AA^1$ and its lift to $X$. Then
$\NN\to \shM_{\AA^1}\subset \O_{\AA^1}$, $1\mapsto t$ is a chart for
the log structure on $\AA^1$. By toric transversality
$\underline\varphi_0(P)$ lies in the big torus of the
$(n-1)$-dimensional toric variety $X_{h(E)}\subset (X_0)_\sing$. The
total space has a singularity of type $A_{e-1}$ along the big torus of
$X_{h(E)}$ where $e$ is the integral length of $h(E)$. Thus via a
toric chart there exist $x,y\in \M_{X,\overline{\underline
\varphi_0(P)}} \subset \O_{X,\overline{\underline\varphi_0(P)}}$ such
that $xy=t^e$ and
\[
S_e=\NN^2\oplus_{\NN} \NN \lra
\O_{X,\overline{\underline\varphi_0(P)}},\quad
\big((a,b),c\big)\longmapsto x^a y^b t^c
\]
is a chart for the log structure on $X$ at $P$. The homomorphism
$\pi^*\shM_{\AA^1}\to \shM_X$ lifts to a map of charts via the
inclusion $\NN\to \NN^2\oplus_{\NN} \NN$ into the second factor. By
the pre-log condition (Definition~\ref{pre-log curves},~(ii)) there
exist generators $z,w\in\O_{C_0,\overline P}$ of the maximal ideal
with
\begin{eqnarray}\label{lift uniformizers}
\underline\varphi_0^*(x)= z^\mu,\quad
\underline\varphi_0^*(y)= w^\mu
\end{eqnarray}
In particular, $zw=0$. Denote by $s_x,s_y,s_t\in\shM_{X_0,\overline
P}$ the restrictions of $x,y,t\in \shM_{X,\overline P}$ to the central
fiber. Then $s_x s_y=s_t^e$ holds in $\shM_{X_0,\overline P}$.
All these choices do not depend on the extension $\shM_{C_0}$ and are
made once and for all.

By (\ref{lift uniformizers}) $\underline\varphi_0^*(s_x)$,
$\underline\varphi_0^*(s_y) \in (\underline\varphi_0^*
\shM_{X_0})_{\overline P}$ map to $z^\mu, w^\mu$ under the structure
homomorphism $(\underline\varphi_0^*\shM_{X_0})_{\overline P}\to
\O_{C_0,\overline P}$, respectively. Thus for any factorization
$\shM_{C_0}\to \O_{C_0}$ of $\underline\varphi_0^*\shM_{X_0}\to
\O_{C_0}$ with the requested properties there exist unique
$s_z,s_w\in \shM_{C_0,\overline P}$ with $\underline\varphi_0^*(s_x)
=s_z^\mu$, $\underline\varphi_0^*(s_y) =s_w^\mu$ and $s_z, s_w$
lifting $z,w$. Then $s_x s_y=s_t^e$ implies $(s_z s_w)^\mu=s_t^e$.
This is equivalent to $s_z s_w=\zeta s_t^{e/\mu}$ for a well-defined
$\mu$-th root of unity $\zeta\in\kk^\times$. This step uses the
integrality of $e/\mu$. Conversely, for any $\zeta\in\kk^\times$ with
$\zeta^\mu=1$ let $\shM_{C_0,\overline P}\to \O_{C_0,\overline P}$ be
the log structure at $\overline P$ defined by the map
\[
S_{e/\mu}\lra \O_{C_0,\overline P},\quad
\big((a,b),c\big)\longmapsto \left\{\begin{array}{ll}
(\zeta^{-1} z)^a w^b,&c=0\\
0,&c\neq0.\end{array}\right.
\]
Denote by $s_{((a,b),c)}\in \shM_{C_0,\overline P}$ the element
belonging to $((a,b),c)\in S_{e/\mu}$. The morphism to $O_0$ given by
mapping $s_t$ to $s_{((0,0),1)}$ is smooth. In particular, its
restriction to $C_0\setminus\{P\}$ is strict and hence equals
$\underline\varphi_0^*\shM_{X_0}|_{C_0\setminus\{P\}}$. Mapping $s_x,
s_y$ to the $\mu$-th powers of $s_z:=\zeta s_{((1,0),0)},
s_w:=s_{((0,1),0)}$, respectively, yields a compatible extension
$(\underline\varphi_0^*\shM_{X_0})_{\overline P}\to
\shM_{C_0,\overline P}$ with associated $\mu$-th root $\zeta$.

Our arguments so far produce a complete list of
$\prod_{E\in\tilde\Gamma^{[1]}\setminus \tilde\Gamma^{[1]}_\infty}
w(E)$ log morphisms $C_0\to X_0$ with the requested properties. It
remains to show that together with the marked points this list
comprises only
\[
w(\Gamma,\mathbf{E}) = \prod_{E\in
\Gamma^{[1]}\setminus \Gamma^{[1]}_\infty} w(E)
\cdot\prod_{i=1}^l w(E_i)
\]
isomorphism classes. First, note that for any of our log structures
$\shM_{C_0}$ on $C_0$ the sheaf $\shM_{C_0}/\O^*_{C_0}$ of finitely
generated monoids is naturally a subsheaf of $\underline\nu_*
\NN_{\underline {\hat C}_0}$ for $\underline\nu: \underline{\hat
C}_0\to \underline C_0$ the normalization. At a nodal point $P\in
C_0$ with chart $S_{e/\mu}\to \O_{C_0,\overline P}$ the image is
generated by $(e/\mu) s_1$, $(e/\mu)s_2$ and $s_1+s_2$ if $s_1$,
$s_2$ are generators of $\underline\nu_* \NN_{\underline{\hat
C}_0,\overline P} = \NN^2$. Away from the nodes it is an isomorphism.
Second, if $\underline\kappa: \underline C_0\to \underline C_0$
commutes with $\underline\varphi_0$ then it maps each irreducible
component and each node to itself. In particular, $\underline\kappa$
induces the identity transformation on $\underline\nu_*
\NN_{\underline{\hat C}_0}$. Thus if $\varphi'_0: C'_0\to X_0$,
$\varphi''_0: C''_0\to X_0$ are two log morphisms from our list and
$\kappa:C''_0\to C'_0$ is an isomorphism with $\varphi''_0=\varphi'_0
\circ\kappa$ then $\kappa$ induces an isomorphism $\kappa^{-1}
\shM_{C'_0}/\O^*_{C_0}\to \shM_{C''_0}/\O^*_{C_0}$. This implies that
$\kappa$ is strict (\cite{F.Kato}, Lemma~3.3).  Therefore it suffices
to investigate the effect of the action of the $\kk$-rational points
of the identity component $\Aut^0(\underline C_0)\subset
\Aut(\underline C_0)$ on our construction of $\shM_{C_0}$ \emph{by
pull-back}. Taking into account the marked points reduces further to
the subgroup $\Aut^0(\underline C_0,\mathbf{x}_0)\subset
\Aut^0(\underline C_0)$.

This action is non-trivial only on the components of $\underline C_0$
having at most two special points, that is, corresponding to
\emph{unmarked} divalent vertices of $\tilde\Gamma$. In particular,
$\Aut^0(\underline C_0) \simeq\GG_m^{\sharp \tilde\Gamma^{[1]}-\sharp
\Gamma^{[1]}-l}$. We are thus reduced to considering the following
situation. Let $E\in\Gamma^{[1]}$, assumed unmarked for the time
being, split into edges $E^1,\ldots,E^r$ in $\tilde\Gamma$ with
$\partial E^i=\{V_{i-1}, V_i\}$ for the bounded edges, and $\partial
E^1=\{V_1\}$ in case $E$ is unbounded (and then
$E^i\in\tilde\Gamma_\infty^{[1]}$ iff $i=1$). Let $P_i\in \underline
C_0$ be the special point labeled by $E^i$. If $E^i$ is bounded then
$P_i$ is a nodal point. In this case the flags $(E^i,V_{i-1})$,
$(E^i,V_i)$ mark the two branches of $\underline C_0$ at $P_i$, for
which we had chosen local coordinates $w_i,z_i \in \O_{C_0,\overline
P_i}$ above. Observe that by their toric construction $w_i,z_i$ extend
to $(\underline C_0)_{V_{i-1}}\setminus \{P_{i-1}\}\cup (\underline
C_0)_{V_i}\setminus \{P_{i+1}\}$. Thus $z_i$, $w_i$ are unique up to
constant and the action of the identity component $\Aut^0(\underline
C_0)\subset \Aut(\underline C_0)$ on $\O_{C_0,P_i}$ becomes linear in
the following sense. Let $\pi_i: \Aut^0(\underline C_0)(\kk) \to
\GG_m(\kk)=\kk^*$ be the projection to the factor acting non-trivially
on the $i$-th irreducible component $(\underline C_0)_{V_i}$. Then for
any $\kappa\in \Aut^0(\underline C_0)(\kk)$
\[
\kappa^* z_r=z_r,\quad
\kappa^* z_i=\pi_i(\kappa)^{-1}\cdot z_i,\quad
\kappa^* w_{i+1}=\pi_i(\kappa)\cdot w_{i+1},\quad
i=1,\ldots, r-1,
\]
and, in the bounded case, $\kappa^* w_1=w_1$. In particular, pulling
back by such $\kappa$ is compatible with our normalizations above iff
$\pi_i(\kappa)^\mu=1$ for all $i$. Thus going inductively from $i=1$
to $r$ we may assume $z_i$, $w_i$ chosen in such a way that the
$\mu$-th root of unity $\zeta_i$ associated to the log structure at
$P_i$ is equal to $1$ for $i=1,\ldots,r-1$. On the other hand, under
pulling back by $\kappa$ the product $\zeta_1\cdot \ldots\cdot
\zeta_r$ remains constant, and hence it is impossible to also
normalize $\zeta_r$. This shows that up to isomorphism there are only
$w(E)$ in the bounded case ($1$ in the unbounded case) rather than
$\prod_{i=1}^r w(E^i)= w(E)^r$ pairwise non-isomorphic choices of
$\shM_{C_0}\to \O_{C_0}$ with the requested properties in a
neighborhood of $\bigcup_{i=1}^r (\underline C_0)_{V_i}\subset
\underline C_0$.

A similar argument holds for the marked edges $E_i$, $i=1,\ldots,l$.
However, since $\Aut^0(\underline C_0,\mathbf{x}_0)$ acts trivially on
the irreducible component $(C_0)_V\subset C_0$ distinguished by the
marked point one has to consider two chains of edges in $\tilde\Gamma$
according to the two connected components of $E\setminus V$. Thus the
total count for $E$ keeps an additional factor of $w(E)$, in agreement
with our definition of total marked weight.
\qed
\medskip

Now that we have a morphism of log spaces $\varphi_0: C_0\to X_0$
over $O_0$ we wish to lift it to the thickenings $O_k:= (\Spec
\kk[t]/(t^{k+1}),\shM_{O_k})$, $k>0$, order by order. For the time
being we forget the marked points and the incidence conditions.
They will be taken care of in a separate discussion. The log
structure on $O_k$ is defined by asking the embedding $O_k\to
\AA^1$ to be strict, that is, by the chart $\NN\to \O_{O_k}$,
$1\mapsto t$. Finding $\varphi_k$ involves both a lift $C_k\to
O_k$ of the domain and an extension of $\varphi_0$ to $C_k$. In
other words we are looking for diagrams of the form
\[\begin{CD}
C_0 @>>> C_k @>>>X\\
@VVV @VVV @VVV\\
O_0@>>> O_k@>>> \AA^1,
\end{CD}\]
with given outer rectangle (formed by $C_0$, $X$, $O_0$ and $\AA^1$)
and lower horizontal sequence. Denote such a diagram by $[C_k/O_k\to
X]$. In the $k$-th step the answer involves the \emph{logarithmic
normal sheaf} $\shN_{\varphi_0}:= \varphi_0^*\Theta_{X/\AA^1}/
\Theta_{C_0/O_0}$. Note that this sheaf is locally free of rank
$n-1$, because pulling-back a differential with logarithmic pole
along a divisor $B$ under a cover with branch locus $B$ preserves the
logarithmic pole and hence $\varphi_0^*\Omega^1_{X/\AA^1} \to
\Omega^1_{C_0/O_0}$ is surjective.

\begin{lemma}\label{deformation theory}
Let $[\varphi_{k-1}: C_{k-1}/O_{k-1}\to X]$ be a lift of
$[\varphi_0:C_0/O_0\to X]$ and assume $C_0$ is rational. Then the set
of isomorphism classes of lifts $[\varphi_k: C_k/O_k\to X]$
restricting to $\varphi_{k-1}$ over $O_{k-1}$ is a torsor under
$H^0(C_0,\shN_{\varphi_0})$.
\end{lemma}
\proof
The obstruction group for the existence of a smooth lift $C_k\to
O_k$ of $C_{k-1}\to O_{k-1}$ is $\Ext^2 (\Omega^1_{C_{k-1}
/O_{k-1}}, \O_{C_0})$, see \cite{K.Kato},~3.14 and \cite{F.Kato},
Proposition~8.6. Since $\Omega^1_{C_{k-1}/O_{k-1}}$ is locally free this
group equals $H^2(C_0,\Theta_{C_0/O_0})$, which vanishes for
dimension reasons. Hence $C_k\to O_k$ exists.

Once a deformation of the domain is given, an extension of the
morphism $\varphi_{k-1}$ exists locally because $X\to \AA^1$ is
log-smooth \cite{K.Kato}, 3.3. We may thus cover $C_k$ by finitely
many open sets $U$ such that a local extension $\varphi_k^U$ exists.
The set of local extensions of $\varphi_{k-1}$ to $C_k$ form a torsor
under the sections of $\varphi_{k-1}^* \Theta_{X/\AA^1}$, see
\cite{K.Kato}, Proposition~3.9. Thus the differences between the
$\varphi_k^U$ on intersections define a \v Cech $1$-cocycle with
values in $\varphi_{k-1}^* \Theta_{X/ \AA^1}$. But $\Theta_{X/\AA^1} =
N\otimes_\ZZ \O_X$ and hence $H^1(C_k,\varphi_{k-1}^* \Theta_{X/
\AA^1})=0$ by rationality of $C_0$. Thus there exists a correction of
the $\varphi_k^U$ by a Cech $0$-cochain to make them patch to a
morphism $\varphi_k$. This shows that the space of liftings
$[\varphi_k: C_k/O_k\to X]$ of $\varphi_{k-1}$ is non-empty. Note that
unlike ordinary deformation theory we are not allowed to take the
trivial deformation because of the requirement of log-smoothness.

To classify the set of isomorphism classes of such liftings knowing
that at least one exists recall first that there is a one-to-one
correspondence between extensions $C_k/O_k$ of $C_{k-1}/O_{k-1}$ and
$\O_{C_{k-1}}$-module extensions $\shE$ of
$\Omega^1_{C_{k-1}/O_{k-1}}$ by $\O_{C_0}$ as follows \cite{K.Kato}.
Given $\shE$ define $\O_{C_k}$ by the fiber product of $\shE\to
\Omega^1_{C_{k-1} /O_{k-1}}$ and of the composition
\[
\O_{C_{k-1}}\stackrel{d}{\lra}
\Omega^1_{\underline C_{k-1}/\underline O_{k-1}}
\lra \Omega^1_{C_{k-1}/O_{k-1}}.
\]
This identifies $\O_{C_k}$ with a subalgebra of the trivial extension
$\O_{C_{k-1}}[\shE]$. Define the extension $\shM_{C_k}$ of the log
structure analogously by fiber product of $\shE\to \Omega^1_{C_{k-1}
/O_{k-1}}$ and of $\dlog: \shM_{C_{k-1}}\to \Omega^1_{C_{k-1}/
O_{k-1}}$. As one easily checks the structure homomorphism
$\shM_{C_{k-1}}\to \O_{C_{k-1}}$ lifts uniquely to $C_k$. Conversely,
any extension $C_k/O_k$ gives rise to two diagrams with exact rows,
one for the extension of scheme structure 
\[\begin{CD}
0@>>> t^k \O_{C_0}@>>> \O_{C_k}@>>> \O_{C_{k-1}}@>>>0\\
@. @| @VV{d}V @VV{d}V\\
0@>>> \O_{C_0}@>>> \Omega^1_{C_k/O_k}\otimes \O_{C_{k-1}}
@>>> \Omega^1_{C_{k-1}/O_{k-1}}@>>>0,
\end{CD}\]
and one for the extension of log structure
\[\begin{CD}
0@>>> t^k \O_{C_0}@>>> \shM_{C_k}@>>> \shM_{C_{k-1}}@>>>0\\
@. @| @VV{\dlog}V @VV{\dlog}V\\
0@>>> \O_{C_0}@>>> \Omega^1_{C_k/O_k}\otimes \O_{C_{k-1}}
@>>> \Omega^1_{C_{k-1}/O_{k-1}}@>>>0.
\end{CD}\]
The two diagrams are identical after restricting the upper horizontal
sequences to invertible elements. In any case, the one-to-one
correspondence thus identifies the sheaf extension $\shE\in
\Ext^1_{C_{k-1}} (\Omega^1_{C_{k-1} /O_{k-1}}, \O_{C_0})$ with the
restriction to $C_{k-1}$ of the sheaf of relative log-differentials
$\Omega^1_{C_k/O_k}$ of the corresponding extension $C_k/O_k$ of
log-spaces.

For the extension of morphism to $X$ observe that an extension of
$\varphi_{k-1}: C_{k-1}\to X$ to $C_k$ over $\AA^1$ gives a
factorization of $\varphi_{k-1}^*\Omega^1_{X/\AA^1}\to
\Omega^1_{C_{k-1} /O_{k-1}}$ over $\Omega^1_{C_k/O_k}
\otimes\O_{C_{k-1}}= \shE$. It is not hard to see that this sets up a
one-to-one correspondence between isomorphism classes of extensions of
$\varphi_{k-1}: C_{k-1}\to X$ to $C_k$ over $\AA^1$ and factorizations
of $\varphi_{k-1}^*\Omega^1_{X/\AA^1}\to \Omega^1_{C_{k-1} /O_{k-1}}$
over $\Omega^1_{C_k/O_k} \otimes\O_{C_{k-1}}= \shE$. In other words,
we need to classify diagrams
\[\begin{CD}
0@>>> \ker \varphi_{k-1}^* @>>> \varphi^*_{k-1}\Omega^1_{X/\AA^1}
@>\varphi_{k-1}^*>> \Omega^1_{C_{k-1}/O_{k-1}} @>>>0\\
@. @VVV @VVV @|\\
0@>>> \O_{C_0} @>>> \shE
@>>> \Omega^1_{C_{k-1}/O_{k-1}} @>>>0
\end{CD}\]
with given upper horizontal sequence, up to isomorphism. Note that the
left square is co-cartesian. A simple exercise in homological algebra
now shows that such diagrams are given uniquely up to isomorphism by
the left vertical homomorphism. The proof is finished by noting
$\Hom(\ker \varphi_{k-1}^*,\O_{C_0})= \Hom(\ker \varphi_0^*,
\O_{C_0}) = H^0(C_0, \shN_{\varphi_0})$ because $\shN_{\varphi_0}$ is
free.
\qed
\medskip

It remains to take care of the incidence conditions given by the
intersections with $Z_i$. To this end we first need to include
deformations of the tuple of marked points $\mathbf x_0
=(x_1,\ldots,x_l)$. In the $k$-th step this merely amounts to go over
from $\Omega^1_{C_k/O_k}$ to the twisted module
$\Omega^1_{C_k/O_k}(\mathbf x)$, which brings in another
$l$-dimensional freedom for the lift. We are then talking about
lifting pairs $[\varphi_{k-1}: C_{k-1}/ O_{k-1}\to X,\mathbf x_{k-1}]$
to $O_k$. In each step we want the tuple of sections $\mathbf x_{k-1}:
\underline O_{k-1}\to \underline X$ to factor over the tuple $\mathbf
Z$ of incidence varieties.

\begin{proposition}\label{deformation with incidence conditions}
Let $[\varphi_{k-1}: C_{k-1}/O_{k-1}\to X,\mathbf x_{k-1}]$ be a lift
of $[\varphi_0:C_0/O_0\to X,\mathbf x_0]$ with $\mathbf x_{k-1}$
factoring over $\mathbf Z$, and assume $C_0$ is rational. Then up to
isomorphism there is a unique lift $[\varphi_k: C_k/O_k\to X,\mathbf
x_k]$ to $O_k$ with $\mathbf x_k$ factoring over $\mathbf Z$.
\end{proposition}

\proof
For a closed point $y\in X$ denote by $T_{X/\AA^1,y}
=\Theta_{X/\AA^1,y}\otimes_{\O_X}\O_{X,y}/\maxid_y$ the fiber of
$\Theta_{X/\AA^1}$ at $y$, which is a $\kk$-vector space of
dimension $n$, and similarly for $\Theta_{Z_i/\AA^1}$ and
$\Theta_{C_0/O_0}$. By Lemma~\ref{deformation theory} it
suffices to prove that the transversality map
\begin{eqnarray}\label{transversality map}
H^0(\shN_{\varphi_0})\lra
\prod_i T_{X/\AA^1,\varphi_0(x_i)}
\big/ \big(T_{Z_i/\AA^1,\varphi_0(x_i)}+
D\varphi_0 (T_{C_0/O_0,x_i})\big)
\end{eqnarray}
is an isomorphism. Now recall that $A_i$ intersects $h(E_i)$ in the
relative interior ($E_i\in \Gamma^{[1]}$ was the $i$-th marked edge).
Hence $A_i\cap h(E_i)$ is a divalent vertex $V$ of $\tilde\Gamma$. In
the discussion in Proposition~\ref{space of lines is torsor} of the
space of lines for this case we saw that $\varphi_0(C_V)$ is then an
orbit closure of $\GG(\ZZ u_i)\subset \GG(N)$ for $u_i\in N$ an
integral generator of $L(h(E_i))\cap N$. This yields a canonical
identification $D\varphi_0 (T_{C_0/O_0,x_i}) = \kk u_i\subset N_\kk:=
N\otimes_\ZZ\kk$. Similarly the actions of $\GG(N)$ and of
$\GG(L(A_i)\cap N)$ identify the other terms on the right-hand side of
(\ref{transversality map}) as follows:
\[
T_{X/\AA^1,\varphi_0(x_i)}= N_\kk,\quad
T_{Z_i/\AA^1,\varphi_0(x_i)}= L(A_i)\otimes_\QQ \kk\subset N_\kk.
\]

To describe $H^0(\shN_{\varphi_0})$ in similar terms observe that the
map $\varphi_0^*\Theta_{X/\AA^1}\to \shN_{\varphi_0}$ induces a
canonical surjection $N\otimes_\ZZ \O_{C_V}\to \shN_{\varphi_0}
\otimes \O_{C_V}$, for any $V\in\tilde\Gamma^{[0]}$. This gives a
canonical surjection $N_\kk= \Gamma (N \otimes_\ZZ \O_{C_V}) \to
\Gamma(\shN_{\varphi_0}\otimes \O_{C_V})$.  To describe the kernel of
this map consider the restriction of the exact sequence
\[
0\lra \Theta_{C_0/O_0}\lra \varphi_0^*\Theta_{X/\AA^1}
\lra \shN_{\varphi_0}\lra 0
\]
to $C_V$. If $V$ is trivalent then $\Theta_{C_0/O_0}|_{C_V}\simeq
\O_{\PP^1}(-1)$. Thus in this case $N_\kk \to\Gamma( \shN_{\varphi_0}
\otimes \O_{C_V})$ is an isomorphism because
$H^i(\Theta_{C_0/O_0}|_{C_V}) = 0$ for $i=0,1$. In the divalent case
$\Theta_{C_0/O_0}|_{C_V}$ is trivial, hence spanned by the vector
field generating the one-parameter subgroup with orbit
$\varphi_0(C_V)$. This shows $\ker \big(N_\kk
\to\Gamma(\shN_{\varphi_0}\otimes \O_{C_V})\big) =\kk u_{(V,E)}$ for
any edge $E$ adjacent to $V$. 

Now if $E\in \tilde \Gamma^{[1]}$ and $\partial E=\{V,V'\}$ then the
images of $h,h'\in N_\kk$ glue to a section of
$H^0(\shN_{\varphi_0})$ over $C_V\cup C_{V'}$ if and only if $h-h'\in
L(h(E))\otimes_\QQ\kk$. In fact, because $\varphi_0$ is transverse to
the toric prime divisor $X_{h(E)}\subset X_V$  (or in $X_{V'}$) the
canonical map $T_{X_{h(E)},\varphi_0(P_E)}\to \shN_{\varphi_0,P_E}$
is an isomorphism. Here $P_E= C_V\cap C_{V'}$ is the nodal point
corresponding to $E$. Thus the claim follows by noting
$T_{X_{h(E)},\varphi_0(P_E)} = (N_\QQ/L(h(E)) ) \otimes_\QQ \kk$. We
conclude that the kernel of
\[
\tilde\Phi: \map(\tilde\Gamma^{[0]},N_\kk)\lra
\prod_{E\in\tilde\Gamma^{[1]}
\setminus \tilde\Gamma_\infty^{[1]}}
N_\kk/\kk u_{(\partial^- E, E)},\quad
h' \longmapsto \big(h'(\partial^+E)-h'(\partial^-E)\big)_E
\]
surjects onto $H^0 (\shN_{\varphi_0})$. Dividing out $\map
(\tilde\Gamma^{[0]}\setminus \Gamma^{[0]}, N_\kk)$ to take into
account the non-uniqueness caused by the divalent vertices, identifies
$H^0 (\shN_{\varphi_0})$ with the kernel of
\[
\Phi: \map(\Gamma^{[0]},N_\kk)\lra
\prod_{E\in\Gamma^{[1]}
\setminus\Gamma_\infty^{[1]}}
N_\kk/\kk u_{(\partial^- E, E)},\quad
h' \longmapsto \big(h'(\partial^+E)-h'(\partial^-E)\big)_E.
\]
The proof is finished by noting that $\Phi$ merged with
(\ref{transversality map}) gives the lattice inclusion (\ref{lattice
map}) from Proposition~\ref{count of maximally degenerate curves},
tensored with $\kk$.
\qed
\medskip

\begin{corollary}\label{formal smoothing}
There exists a unique $l$-marked stable map $(\underline
\varphi_\infty: \underline C_\infty\to \underline X,\mathbf
x_\infty)$ over $\underline O_\infty=\Spec \kk\lfor t\rfor$ incident
to $Z_1,\ldots,Z_l$ and such that the induced log morphism of the
central fibers agrees with $[\varphi_0:C_0/O_0\to X_0,\mathbf x_0]$.
\end{corollary}

\proof
Forgetting the log structure from $[\varphi_k:C_k/O_k\to X_k,\mathbf
x_k]$ in the proposition gives an inverse system of maps from
$\underline O_k$ to the stack of stable maps to $\underline X$. This
stack has been constructed in \cite{behrend manin}. Going
over to the limit provides a map from $\underline O_\infty$ to this
stack, hence a stable map over $\underline O_\infty$.

Conversely, let $(\underline \varphi_\infty: \underline C_\infty\to
\underline X,\mathbf x_\infty)$ be a stable map over $\underline
O_\infty$. On the involved spaces $\underline C_\infty,\underline
O_\infty, X, \AA^1$ take the log structures $\shM_{C_\infty}$,
$\shM_{O_\infty}$, $\shM_X$, $\shM_{\AA^1}$ defined by the divisors
$C_0\subset C_\infty$, $X_0\subset X$, and the closed points $O_0\in
O_\infty$, $0\in \AA^1$, respectively. Because $\underline
\varphi_\infty$ maps the pair $(\underline C_\infty,\underline C_0)$
to $(\underline X,\underline X_0)$ it induces uniquely the commutative
diagram of log spaces
\[\begin{CD}
C_\infty@>>> X\\
@VVV @VVV\\
O_\infty@>>> \AA^1.
\end{CD}\]
From \cite{F.Kato2} it is clear that $C_\infty\to O_\infty$ is log
smooth. By construction the restriction $\varphi_0: C_0\to X_0$ to the
central fiber also fulfills the other properties listed in
Proposition~\ref{log structures on degenerate curves}.
\qed

\begin{remark}
The same arguments work in the higher genus case for trivalent
tropical curves of the correct genus and with deformation space of
the expected dimension $e+(n-3)(1-g)$ (not superabundant). In the
superabundant case one obtains an obstruction bundle that is an input
to the virtual intersection formalism.
\end{remark}

\section{The Main Theorem}
\label{section main theorem}

Throughout this section fix the following data. (1) A complete fan
$\Sigma$ in $N_\QQ$. (2) A map $\Delta:N\setminus\{0\}\to \NN$ with
finite support and such that $\Delta(v)\neq 0$ $\Rightarrow$
$\QQ_{\ge0} v\in\Sigma^{[1]}$. (3) A tuple $\mathbf
L=(L_1,\ldots,L_l)$ of linear subspaces of $N_\QQ$ with $\codim
L_i=d_i+1$ and such that $\sum_{i=1}^l d_i=e+n-3$, $n=\dim_\QQ
N_\QQ$. To this data we are going to associate two numbers, one
counting algebraic curves on a toric variety, the other counting
tropical curves in $N_\QQ$. The result is that both these numbers are
well-defined and agree.

For the tropical count choose $\mathbf A=(A_1,\ldots,A_l)$ with
$A_i\in N_\QQ/ L_i$ in general position for $\Delta$, as provided by
Proposition~\ref{tropical transversality}. Then any tropical curve of
degree $\Delta$ matching $\mathbf A$ fulfills the
conditions~(i)--(iii) of Definition~\ref{non-general tropically}.
Denote by $\T_{0,l,\Delta}$ the set of isomorphism classes of
$l$-marked tropical curves of genus $0$ and degree $\Delta$. Then the
subset
\[
\T_{0,l,\Delta}(\mathbf A) =\big\{(\Gamma,\mathbf E,h)\in
\T_{0,l,\Delta} \,\big|\, h(E_i)\cap A_i\neq\emptyset\big\}
\]
of curves matching $\mathbf A$ is finite. In fact, by
Proposition~\ref{finiteness of types} there are only finitely many
types of tropical curves of fixed genus and degree, while
Proposition~\ref{tropical transversality} says that for general
$\mathbf A$ any two tropical curves of the same type and matching
$\mathbf A$ are isomorphic. For $(\Gamma,\mathbf E,h)\in
\T_{0,l,\Delta}(\mathbf A)$ let $\mathfrak D(\Gamma,\mathbf
E,h,\mathbf A)$ be the lattice index defined in Proposition~\ref{count
of maximally degenerate curves} and
\[
\tilde {\mathfrak D}(\Gamma,\mathbf E,h,\mathbf A):=
\mathfrak D(\Gamma,\mathbf E,h,\mathbf A)
\cdot \prod_{i=1}^l \delta_i,
\]
where $\delta_i$ is the index of $\ZZ u_{(\partial^- E_i,E_i)} +
L(A_i)\cap N$ in $(\QQ u_{(\partial^- E_i,E_i)}+L(A_i))\cap N$ (see
Remark~\ref{maximal degenerate curves remark}). With these conventions
we can define the count on the tropical side as follows.

\begin{definition}
Assume the affine constraint $\mathbf{A}$ is general for $\Delta$. Then
the \emph{number of tropical curves of degree $\Delta$ and genus $0$
matching $\mathbf A$} is
\begin{center}\fbox{
$\displaystyle N_{0,\Delta}^\text{trop}(\mathbf A):=
\sum_{(\Gamma,\mathbf E,h)\in \T_{0,l,\Delta}(\mathbf A)}
w(\Gamma,\mathbf{E})\cdot \tilde{\mathfrak D} (\Gamma,\mathbf
E,h,\mathbf A).$}
\end{center}
\end{definition}

For the count of algebraic curves consider the toric variety
$X(\Sigma)$ defined by $\Sigma$. The toric prime divisors on
$X(\Sigma)$ are denoted $D_v$ with $v\in N$ primitive generators of
the rays of $\Sigma$. For a torically transverse
(Definition~\ref{torically transverse}) stable map $\varphi:C\to
X(\Sigma)$ define a map $\Delta(\varphi):N\setminus\{0\}\to \NN$ as
follows. For primitive $v\in N$ and $\lambda\in \NN$ map
$\lambda\cdot v$ to $0$ if $\QQ_{\ge 0} v\not\in \Sigma^{[1]}$, and
to the number of points of multiplicity $\lambda$ in $\varphi^* D_v$
otherwise. So $\Delta(\varphi)$ counts the number of intersection
points with the toric prime divisors of any given multiplicity.
Clearly $\Delta(\varphi)$ has finite support. Now define
$M_{0,l,\Delta}$ as the set of isomorphism classes of $l$-marked
stable maps $(\varphi:C\to X(\Sigma),\mathbf x)$ of genus $0$
\emph{with $\varphi$  torically transverse} and such that
$\Delta(\varphi) =\Delta$. Note that this space of stable maps
disregards any curves intersecting the toric strata of codimension
$\ge 2$ and hence generally is non-complete.

To bring in the incidence conditions let $P_1,\ldots,P_l\in
X(\Sigma)$ be closed points in the big torus. Define $Z_i=Z_i(P_i):=
\overline{\GG(L_i\cap N). P_i}$ the closure of the orbit through
$P_i$ for the subgroup defined by $L_i\subset N_\QQ$ and denote
$\mathbf Z=\mathbf Z(\mathbf P)=(Z_1,\ldots,Z_l)$. The set of
isomorphism classes of constrained curves is
\[
M_{0,l,\Delta}(\mathbf Z):=\big\{(\varphi:C\to X(\Sigma),\mathbf x) \in
M_{0,l,\Delta} \,\big|\, \varphi(x_i)\in \Int(Z_i)\big\}.
\]
We will see in the Main Theorem that for general $P_1,\ldots, P_l$
this set is finite. Because the $M_{0,l,\Delta}(\mathbf Z)$ fit into
one algebraic family for varying incidence conditions $P_1,\ldots,
P_l$ its cardinality is then the same for any two general choices of
points.

\begin{definition}
The \emph{number of rational curves in $X(\Sigma)$ of degree $\Delta$
and constrained by $\mathbf L$} is
\begin{center}\fbox{
$\displaystyle N_{0,\Delta}^\text{alg}(\mathbf L):=
\sharp\, M_{0,l,\Delta}\big(\mathbf Z(\mathbf P)\big)$}
\end{center}
for general $\mathbf P=(P_1,\ldots,P_l)$.
\end{definition}

No signs or multiplicities enter in this definition, so
$N_{0,\Delta}^\text{alg}(\mathbf L)$ is a true enumerative count.

\begin{theorem}\label{main theorem}
There exist $P_1,\ldots,P_l\in X(\Sigma)$ such that
$M_{0,l,\Delta}\big(\mathbf Z(\mathbf P)\big)$ is finite, and hence
$N_{0,\Delta}^\text{alg}(\mathbf L)$ is well defined. Moreover, it
holds
\begin{center}\fbox{
$\displaystyle N_{0,\Delta}^\text{alg}(\mathbf L)=
N_{0,\Delta}^\text{trop}(\mathbf A).$}
\end{center}
\end{theorem}

\begin{corollary}\label{tropical count depends only on L}
$N_{0,\Delta}^\text{trop}(\mathbf A)$ depends only on $\Delta$ and
$\mathbf L$.
\end{corollary}

\begin{remark}
It is worthwhile to remark that $N_{0,\Delta}^\text{alg}(\mathbf L)$
is \emph{not} in general a Gromov-Witten invariant in the traditional
sense, and in particular, it may not be a symplectic invariant. This
is due to the fact that our count disregards stable maps with
components mapping to a toric divisor. For example, the results on
$\PP^1\times\PP^1$ differ from the results on the symplectomorphic
even degree Hirzebruch surfaces $\FF_{2k}$.

On the other hand, we believe that our numbers should have an
interpretation as \emph{relative Gromov-Witten invariant} for the
pair $(X(\Sigma), D(\Sigma))$ where $D(\Sigma)= X(\Sigma)\setminus
\GG(N)$ is the union of the toric prime divisors. Unfortunately a
theory of relative Gromov-Witten invariants so far exists only for
smooth divisors and hence the verification of this claim remains
elusive at this point for $n\ge 2$. Nevertheless, our statement is
true provided a theory of relative Gromov-Witten invariants for
pairs $(X,D)$ has the following expected property.
\smallskip

\begin{center}\begin{minipage}{14cm}
\emph{If every stable map of the considered type, fulfilling the
incidence conditions and without components mapping to $D$, intersects
$D$ in the smooth locus and is infinitesimally rigid (unobstructed
deformations) then the relative Gromov-Witten invariant coincides
with the number of such curves.}
\end{minipage}\end{center}
\smallskip

\noindent
In fact, we will see in the proof of the theorem that for general
$P_1,\ldots, P_l$ the assumptions in this statement are fulfilled in
our case.
\qed
\end{remark}

\noindent
\emph{Proof of the theorem.}
We will show that both numbers agree with the number of isomorphism
classes of maximally degenerate curves together with log smooth
structures, for a sufficiently fine toric degeneration of
$X(\Sigma)$.\\[1ex]
{\sl 1) Construction of degeneration.} To define an appropriate toric
degeneration let $S$ be the intersection of $A_1\cup\ldots\cup A_l$
with the union of the images of all tropical curves in
$\T_{0,l,\Delta}(\mathbf A)$. This set is finite for
$\T_{0,l,\Delta}(\mathbf A)$ is finite and the intersections with the
$A_i$ are transverse. Apply Proposition~\ref{adapted polyhedral
decomposition} to the disjoint union of the finitely many tropical
curves in $\T_{0,l,\Delta}(\mathbf A)$, with the requirement that $S$
be contained in the set of vertices. This gives a polyhedral
decomposition $\P$ of $N_\QQ$ with asymptotic fan $\Sigma=\Sigma_\P$,
with $S\subset\P^{[0]}$ and such that for all $(\Gamma,\mathbf
E,h)\in \T_{0,l,\Delta}(\mathbf A)$
\[
h(\Gamma^{[\mu]})\subset \bigcup_{\Xi\in\P^{[\mu]}} \Xi,\quad
\mu=0,1,\quad\text{ and }\quad h(\Gamma)\cap A_i\subset
\P^{[0]},\quad j=1,\ldots,l.
\]
After rescaling we may assume $\P$ to be integral, that is
$\P^{[0]}\subset N$. Then the image of each bounded edge of a
tropical curve has an integral length, the number of integral points
on this image minus one. After another rescaling we may assume that
this integral length is a multiple of the weight of this edge, for
any bounded edge of any tropical curve in $\T_{0,l,\Delta}(\mathbf
A)$. This gives our polyhedral decomposition $\P$ and an associated
toric degeneration $\pi:X=X(\widetilde\Sigma_\P)\to \AA^1$ with
general fiber isomorphic to $X(\Sigma)$ and reduced special fiber
$X_0$.

As for the incidence conditions $Z_i\subset X(\Sigma)$ consider the
subgroup $\tilde\GG_i:=\GG\big(LC(A_i)\cap (N\times\ZZ)\big)\subset
\GG(N\times\ZZ)$ associated to $A_i\subset N_\QQ$. Identify
$X(\Sigma)$ with any general fiber $\pi^{-1} (Q)$ and define $\tilde
Z_i$ as the closure of the $\tilde \GG_i$-orbit through $P_i$. Then
$\tilde Z_i\cap X(\Sigma)= Z_i$, so $\tilde Z_i$ is a degeneration of
$Z_i$. If $v\in A_i\cap \P^{[0]}$ then the intersection of $\tilde
Z_i$ with the irreducible component $X_v$ is the closure of a
$\GG(L_i\cap N)$-orbit of maximal dimension, that is, intersecting the
big torus (Proposition~\ref{orbit closures for affine constraints}).
Along the big torus this intersection is even transverse
(Corollary~\ref{transversality of incidence conditions}).\\[1ex]
{\sl 2) Interpretation of $N_{0,\Delta}^\text{trop} (\mathbf A)$.} The
number of isomorphism classes of maximally degenerate curves in $X_0$
intersecting $\tilde Z_i\cap X_0$ and with dual intersection graph
$(\Gamma,h)\in \T_{0,l,\Delta}(\mathbf A)$ equals ${\mathfrak
D}(\Gamma,\mathbf E,h,\mathbf A)$ (Proposition~\ref{count of maximally
degenerate curves}). According to Remark~\ref{maximal degenerate
curves remark} the image of each such stable map $\varphi_0: C_0\to
X_0$ intersects $Z_i$ in $\delta_i$ points. Recall that the
irreducible components of $C_0$ corresponding to the point of
intersection $v\in h(E_i)\cap A_i$ is a copy of $\PP^1$ mapping to a
line in $X_v$ by a $w(E_i)$-fold covering, unbranched over
$\Int(X_v)$. Thus each of the $\delta_i$ intersection points has
$w(E_i)$ preimages on $C_0$. Any of these is a possible choice for the
$i$-th marked point ${x_0}_i\in C_0$. However, by the same token,
choices lying in the same $\varphi_0$-fiber lead to isomorphic stable
maps. Thus up to isomorphism there are $\tilde {\mathfrak
D}(\Gamma,\mathbf{E}, h,\mathbf A) = {\mathfrak D}(\Gamma,\mathbf{E},
h,\mathbf A)\cdot \prod_i\delta_i$ \emph{marked} stable maps
$(C_0,\mathbf x,\varphi_0)$ with marked dual intersection graph
$(\Gamma,\mathbf{E},h)$. Choose one of them.

According to Proposition~\ref{log structures on degenerate curves}
there are $w(\Gamma,\mathbf E)$ pairwise non-isomorphic ways to lift
$\varphi_0$ to a morphism of log smooth spaces $(C_0,\M_{C_0}) \to
(X_0,\M_{X_0})$ relative the standard log point $(\Spec \kk,
\NN\times\kk^\times)$. Hence by definition $N_{0,\Delta}^\text{trop}
(\mathbf A)$ agrees with the total number of isomorphism classes of
such maximally degenerate stable maps with $l$ marked points mapping
to $\tilde Z_1,\ldots, \tilde Z_l$, and together with log structures.
\\[1ex]
{\sl 3) Interpretation of $N_{0,\Delta}^\text{alg} (\mathbf L)$.} The
results from Section~\ref{section deformation theory} provide the
link to curves on $X(\Sigma)$. Corollary~\ref{formal smoothing} says
that each log morphism $(C_0,\M_{C_0}) \to (X,\M_{X})$ is induced by
a unique family of stable maps $(C_\infty\to X,\mathbf x_\infty)$
over $\Spec \kk\lfor t\rfor$, both incident to $\tilde
Z_1,\ldots,\tilde Z_l$. For the relevant stack of $l$-marked stable
maps to $X$ over $\AA^1$ with incidences this means that there is an
injection from the set of log morphisms $(C_0,\M_{C_0}) \to
(X,\M_{X},\mathbf x_0)$ considered in Section~\ref{section
deformation theory} to the set of prime components of this
Deligne-Mumford stack at maximally degenerate curves mapping
dominantly to $\AA^1$. We claim that every such prime component
arises in this way.

To this end assume that $R$ is a discrete valuation ring with residue
field $\kk$ and $\Spec R\to \AA^1$ is a morphism mapping the closed
point to $0\in\AA^1$. This map identifies the completion of $R$ with
$\kk\lfor t\rfor$. Let $\varphi^*: C^*\to X\setminus X_0$ together
with $l$ sections $x^*_1,\ldots,x^*_l$ of marked points be a stable
map defined over the quotient field $\Spec K\subset \Spec R$.
Corollary~\ref{stable reduction} says that possibly after a ramified
base change and toric modification $\tilde X\to X$ with centers on the
central fiber, $\varphi^*$ together with the $l$ sections extends to a
stable map $(\varphi,\mathbf x)$ defined over $\Spec R$, such that for
every irreducible component $X_v\subset X_0$ the projection
$C\times_{X_0}X_v \to X_v$ is a torically transverse stable map. We
will see towards the end of the proof that this base change and
modification are indeed unnecessary. For the time being let us just
make this transformation but keep the notations $\P$, $X$, etc.\ as
above. Our goal is to show that the completion of $\varphi$ at the
closed point of $\Spec R$ is already in the list of stable maps
constructed by Corollary~\ref{formal smoothing} above.

Let us first check that $\varphi_0$ is indeed a pre-log curve
(Definition~\ref{pre-log curves}). In fact, if $P\in C_0$ is a closed
point mapping to the singular locus of $X_0$ consider the homomorphism
of complete local $\kk\lfor t\rfor$-algebras $\psi:\hat
\O_{X,\varphi_0(P)}\to \hat \O_{C,P}$. The following possibilities
arise: (1) $C_0$ is smooth at $P$; (2) $C_0$ has a node at $P$ with
both branches mapping to the same prime component of
$(X_0,\varphi_0(P))$; (3) $C_0$ has a node at $P$, but the two
branches map to different prime components of $(X_0,\varphi_0(P))$. In
the first case $C$ is also smooth at $P$, and by the maximal degeneracy
condition $\psi$ has the following form:
\[
\kk\lfor x,y,t,u_1,\ldots,u_{n-1}\rfor/(xy-t^e)\lra \kk\lfor z,t\rfor,
\quad \psi(x)= z^\alpha g+th,\ \psi(y)= tk,
\]
with $g$ a unit. But then $t^e=\psi(xy)= (z^\alpha g+th)tk$, and comparing
monomials shows that this case can not occur. A
similar argument shows the impossibility of (2): By Lemma~\ref{Hom of
A_k singularities},~(1) below, $\hat \O_{C,P}\simeq \kk\lfor
z,w,t\rfor/(zw-\lambda t^b)$ with $\lambda\in\{0,1\}$ and $b>0$. Then
$\psi$ takes the form
\[
\kk\lfor x,y,t,u_1,\ldots,u_{n-1}\rfor/(xy-t^e)\lra
\kk\lfor z,w,t\rfor/(zw-\lambda t^b),\quad
\psi(x)= z^\alpha g+tk,\ \psi(y)= z^\beta h+tl,
\]
with $g,h$ units. Again this leads to a contradiction with
$\psi(xy)=t^e$. We are thus left with~(3). In this case
Lemma~\ref{Hom of A_k singularities},(2) below shows that the
intersection numbers of the two branches of $\varphi_0$ with the
singular locus of $X_0$ coincide. Hence Condition~(ii) in
Definition~\ref{pre-log curves} is fulfilled.

Thus by Construction~\ref{open dual intersection graph} the map from
the dual intersection complex of $C_0$ to the $1$-skeleton of $\P$
defines a tropical curve $h:\Gamma\to N_\RR$ of genus $0$ and degree
$\Delta$, and a corresponding subdivision $\tilde\Gamma$ with
$\tilde\Gamma^{[0]}= h^{-1}(\P^{[0]})$. It also comes with incidences
with $A_i$ as follows. Let $V_i\in \tilde\Gamma^{[0]}$ be the vertex
corresponding to the irreducible component of $C_0$ containing the
$i$-th marked point. Because $\varphi\circ x_i:\Spec R\to X$ factors
over the $i$-th incidence condition $\tilde Z_i$ it follows
$\Int(X_{h(V_i)})\cap \tilde Z_i\neq\emptyset$. Recalling that $\tilde Z_i$
is an orbit closure for $\GG(LC(A_i)\cap(N\times \ZZ))$
Proposition~\ref{orbit closures for affine constraints} implies
$h(V_i)\in A_i$.

We have now shown that $(\Gamma,h,\mathbf E)\in \T_{0,l,\Delta}(A)$.
In other words, it is one of the tropical curves already considered
above, and hence it is general for $\mathbf A$ in the sense of
Definition~\ref{non-general tropically}. In particular, all vertices
of $\Gamma$ are at most trivalent, the $V_i$ are divalent and the map
$\Gamma^{[1]}\to \P^{[1]}$ induced by $h$ is injective. For $n>2$ the
map $h$ itself is injective. For the following discussion let us make
this assumption and leave the straightforward modifications for the
case when $h$ is only an immersion to the interested reader. Then for
every irreducible component $X_v\subset X_0$ there is at most one
irreducible component $C_V\subset C_0$ mapping to $X_v$, namely if
$v=h(V)$ for $V\in \tilde\Gamma^{[0]}$, and $\varphi_0|_{C_V}$
intersects the toric boundary of $X_v$ in at most three points. Hence
$\varphi_0|_{C_V}$ is a line for every $V\in \Gamma^{[0]}$ and
$\varphi_0$ is indeed maximally degenerate. Because it is of the
correct genus and degree, and because it is incident to $\tilde
Z_i\cap X_0$ for every $i$, it is one of the maximally degenerate
curves constructed for $(\Gamma,h,\mathbf E)\in
\T_{0,l,\Delta}(\mathbf A)$ above (for the refined degeneration
$\tilde X$). The formal completion at the maximal ideal of $R$ thus
defines a family of stable maps $(C_\infty\to X,\mathbf x_\infty)$
over $\kk\lfor t\rfor$ isomorphic to one of the families constructed
by log deformation theory before.
\qed
\medskip

In the proof we used the following algebraic results.

\begin{lemma}\label{Hom of A_k singularities}
1) A $\kk\lfor t\rfor$-algebra of the form $\kk\lfor z,w,t\rfor
/(zw-t^a f)$ with $a>0$ is isomorphic to $\kk\lfor
z',w',t\rfor/(z'w'-\lambda t^b)$ with $\lambda\in\{0,1\}$, $z'=zg$,
$w'=wh$ for units $g,h$ and some $b\ge a$.\\[1ex]
2) If $\phi: \kk\lfor x,y,t,u_1,\ldots,u_{n-1}\rfor /(xy-t^e )\to
\kk\lfor z,w,t\rfor /(zw-t^a f)$ is a $\kk\lfor t\rfor$-algebra
homomorphism with $\phi(x)= z^\alpha\cdot g$, $\phi(y)=w^\beta\cdot
h$, $g,h$ units, and $\alpha,\beta,a,e>0$ then $\alpha=\beta$.
\end{lemma}
\proof
1) This is well-known (and elementary to prove).\\[1ex]
2) The composition with the inclusion
\[
\kk\lfor x,y,t\rfor/(xy-t^e )\lra
\kk\lfor x,y,t,u_1,\ldots,u_{n-1}\rfor/(xy-t^e )
\]
reduces to the case $n=1$. We may also assume $\alpha\le \beta$ and,
by (1), $f=1$. In fact, if $f=0$ then $xy$ maps to zero, but $xy=t^e$,
so this case can not arise. Now $t^e=xy$ maps to a unit times
$z^\alpha w^\beta= (zw)^\alpha w^{\beta-\alpha} =t^{a\alpha}
w^{\beta-\alpha}$. Both $t^e$ and $t^{a\alpha} w^{\beta-\alpha}$ are
part of the $\kk$-vector space basis of $\kk\lfor z,w,t\rfor/
(zw-t^a)$ consisting of all monomials without $zw$-factor,
and hence $a\alpha= e$ and $\beta=\alpha$.
\qed

\begin{remark}
The proof of the theorem gives a more refined correspondence between
the tropical and geometric counts, once a sufficiently fine
degeneration $X\to \AA^1$ has been chosen. First, a choice of tropical
curve matching $\mathbf{A}$ gives combinatorial information for the
degeneration of stable maps, via intersection information on $X_0$.
Then ${\mathfrak D}(\Gamma,\mathbf E,h,\mathbf A)$ is the number of
different degenerate \emph{unmarked} stable maps. Next, there are
$\delta_i$ points of intersection of the image with $Z_i$ compatible
with the marked tropical curve. Hence $\tilde {\mathfrak
D}(\Gamma,\mathbf E,h,\mathbf A) = {\mathfrak D}(\Gamma,\mathbf
E,h,\mathbf A) \cdot \prod_i \delta_i$ is the number of different
central fibers as \emph{marked} stable maps. There are  $w(\Gamma)$
ways of endowing such an unmarked stable map to $X_0$ with a log
structure coming from a degeneration; taking into account the markings
reduced the automorphism group to give $w(\Gamma,\mathbf{E})=
w(\Gamma)\cdot \prod_i w(E_i)$ possibilities. Finally, there is a
one-to-one correspondence between such marked stable maps with log
structures and degenerations of marked stable maps of the considered
type.
\end{remark}
\medskip

In two dimensions Mikhalkin gave a different definition of $N_{0,
\Delta}^\text{trop}$ (\cite{mikhalkin}, Definition~4.16), while his
definition on the geometric side (\cite{mikhalkin}, Definition~5.1)
agrees with ours. The definitions on the tropical sides must therefore
also coincide. We end this section with an independent proof of this
statement. Note that point constraints determine the corresponding
marking of a matching tropical curve uniquely, so in dimension two
markings are redundant information.

Recall that Mikhalkin defines the \emph{multiplicity} of a tropical
curve $(\Gamma,h)$ as the product over the multiplicities
$\mult(\Gamma,h,V)$ of the trivalent vertices. If $E_1$, $E_2$ are two
different edges emanating from a trivalent vertex $V$ then
\[
\mult(\Gamma,h,V)= w(E_1) w(E_2)
\cdot\big|\det(u_{(V,E_1)},u_{(V,E_2)})\big|.
\]
By the balancing conditon this number does not depend on the choices
of $E_1$ and $E_2$. The claimed equivalence of the definition in
\cite{mikhalkin} with ours follows readily from the following result.

\begin{proposition}
Suppose that $(\Gamma,h,\mathbf{E})$ is a genus zero tropical curve
that is general for a tuple $\mathbf{P}=(P_1,\ldots,P_l)$ of points
in $N_\QQ\simeq \QQ^2$. Then
\[
w(\Gamma)\cdot {\mathfrak D}(\Gamma,\mathbf{E},h,\mathbf{P})
=\prod_{V\in\Gamma^{[0]}} \mult(\Gamma,h,V).
\]
\end{proposition}

\proof
The proof is by induction on the number of vertices. If there is only
one vertex there are two point constraints. Letting $E_i$ be the edge
with $P_i \in h(E_i)$ for the computation of $\mult(V)$ shows that
both sides agree trivially.

In the general case with $l$ marked points there are $l+1$ unbounded
edges, $l-1$ vertices and $l-2$ bounded edges. In particular, there
must be one unmarked, unbounded edge $E$. Removing $E$ splits
$\Gamma$ into two connected components, say $\Gamma_1$ and
$\Gamma_2$. Let us first treat the case that both $\Gamma_i$ have
vertices. In this case the restrictions of $(\Gamma,h,\mathbf{E})$ to
the connected components define tropical curves
$(\Gamma_i,h_i,\mathbf{E}_i)$ themselves, and there is a splitting of
$\mathbf{P}$ into two tuples of points $\mathbf{P}_1$, $\mathbf{P}_2$
such that $(\Gamma_i,h_i,\mathbf{E}_i)$ is general for
$\mathbf{P}_i$. Otherwise one $(\Gamma_i,h_i,\mathbf{E}_i)$ moves in
a one-parameter family with each member matching $\mathbf{P}_i$,
contradicting the generality of $(\Gamma,h,\mathbf{E})$.

Let $b_1,\ldots, b_{l-2}\in N$ be primitive generators of the lines
spanned by the images of the bounded edges of $\Gamma$. We order the
$b_i$ in such a way that the first $s-2$ are the images of the bounded
edges of $\Gamma_1$, $b_{s-1}$ and $b_s$ belong to the unbounded edges
$E'$ and $E''$ of $\Gamma_1$ and $\Gamma_2$ with closure in $\Gamma$
intersecting $E$, and the last $l-2-s$ are the images of the bounded
edges of $\Gamma_2$. Similarly, let $u_1,\ldots,u_l$ denote generators
of the lines spanned by $h(E_1),\ldots,h(E_l)$ with the first $s$
coming from $\Gamma_1$ and the rest coming from $\Gamma_2$. Now
${\mathfrak D}_1={\mathfrak
D}(\Gamma_1,\mathbf{E}_1,h_1,\mathbf{P}_1)$ is the order of the
cokernel of a map
\[
\Phi': \map(\Gamma_1^{[0]},N)\lra \prod_{i=1}^{s-2} N/\ZZ b_i\times
\prod_{j=1}^s N/\ZZ u_j,
\]
and ${\mathfrak D}_2={\mathfrak D}(\Gamma_2,\mathbf{E}_2,h_2,\mathbf{P}_2)$ is
the order of the cokernel of a map
\[
\Phi'': \map(\Gamma_2^{[0]},N)\lra \prod_{i=s+1}^{l-2} N/\ZZ b_i\times
\prod_{j=s+1}^l N/\ZZ u_j.
\]
On the other hand, ${\mathfrak D}(\Gamma,\mathbf{E},h,\mathbf{P})$ is then the
order of the cokernel of
\begin{eqnarray*}
\Phi: \map(\Gamma_1^{[0]},N)\times \map(\Gamma_2^{[0]},N)\times N
&\lra&B'\times B''\times N/\ZZ b_{s-1}\times N/\ZZ b_s\\
(h',h'',Q)&\longmapsto&
\big(\Phi'(h'),\Phi''(h''),
\overline{h'(V')- Q},\overline{h''(V'')-Q} \big).
\end{eqnarray*}
Here $V'\in\Gamma_1^{[0]}$, $V''\in \Gamma_2^{[1]}$ are the unique
vertices adjacent to $E'$ and $E''$, respectively, and $B'$, $B''$
denote the free abelian groups on the right-hand sides of the
previous two displayed equations. The overlining indicates taking
equivalence classes. It is then not hard to see that
\[
\big|\coker(\Phi)\big|=\big|\coker(\Phi')\big|\cdot
\big|\coker(\Phi'')\big| \cdot\big|\coker\big(N\to N/\ZZ b_{s-1}\times
N/\ZZ b_s \big)\big|.
\]
The last term on the right-hand side gives $\det(b_{s-1},b_s)$. Now
it holds $w(\Gamma)=w(\Gamma_1)\cdot w(\Gamma_2)\cdot w(E')\cdot
w(E'')$ and hence, taking into account the induction hypothesis
applied to $(\Gamma_i,h_i,\mathbf{E}_i)$,
\begin{eqnarray*}
w(\Gamma)\cdot {\mathfrak D}(\Gamma,\mathbf{E},h,\mathbf{P})
&=&w(\Gamma_1) {\mathfrak D}_1\cdot
w(\Gamma_2) {\mathfrak D}_2\cdot w(E')w(E'')\det(b_{s-1},b_s)\\
&=& \prod_{V\in\Gamma_1^{[0]}} \mult(\Gamma_1,h_1,V)\cdot
\prod_{V\in\Gamma_2^{[0]}} \mult(\Gamma_2,h_2,V)
\cdot \mult(\Gamma,h,\hat V).
\end{eqnarray*}
Here $\hat V\in \Gamma^{[0]}$ is the vertex separating $\Gamma_1$ and
$\Gamma_2$ in $\Gamma$. The expression in the last line equals
$\prod_{V\in\Gamma^{[0]}} \mult(\Gamma,h,V)$ as claimed.

In the case where one of $\Gamma_i$ has no vertices, say $\Gamma_2$,
by generality $\Gamma_2$ must contain a marked point, say the last
one. If $\Phi'$, $B'$, $V'\in\Gamma_1$ are defined as above the
definition for $\Phi$ now reads
\begin{eqnarray*}
\Phi: \map(\Gamma_1^{[0]},N)\times N
&\lra&B'\times N/\ZZ b_{l-2}\times N/\ZZ u_l\\
(h',Q)&\longmapsto&
\big(\Phi'(h'),\overline{h'(V')-Q},\overline{Q} \big).
\end{eqnarray*}
Again this relates $\big|\coker(\Phi)\big|$ and
$\big|\coker(\Phi')\big|$ by the determinant $\det(b_{l-2},u_l)$ of
the two generators of the lines emanating from the image of the
removed vertex. The computation is finished by the inductive
assumption as before. Note that while the edge coming from $\Gamma_2$
is now unbounded, its weight still contributes to our definition of
$N_{0,\Delta}^\text{trop}$ because it contains a marked point.
\qed

\section{Examples}
The enumeration of tropical curves of given degree and fixed genus is
a finite problem. However, the algorithmic problems involved in doing
this efficiently are not well-studied except in dimension two, where
Mikhalkin gave a beautiful algorithm via counting of lattice paths
\cite{mikhalkin}.

In higher dimensions the combinatorics gets out of hands quickly, so
that with a naive approach one is limited to tropical curves with
few ends or to very special geometry. The following observation often
drastically reduces the number of cases to be considered in higher
dimensions. Define the \emph{distance} of two edges in a graph as the
minimal number of other edges that one has to follow in any connecting
path. So edges with a common vertex have distance $0$.

\begin{proposition}\label{distance of incidence conditions}
Let an affine constraint $\mathbf A=(A_1,\ldots,A_l)$ of codimension
$(d_1,\ldots,d_l)$ be general for a type $\Delta$. Then for any
tropical curve $(\Gamma,\mathbf E, h)$ matching $\mathbf A$ the
distance between $E_i$ and $E_j$ is at least $d_i+d_j-n$, for any
$i\neq j$.
\end{proposition}
\proof
In a path $E_i=E(0), E(1),\ldots, E(\delta), E(\delta+1)=E_j$
connecting $E_i$ and $E_j$ and passing through vertices $V_\nu\in
E(\nu)\cap E(\nu+1)$ the $i$-th incidence condition restricts
$h(V_0)$ to an $(n-d_i)$-dimensional subspace. Inductively,
$h(V_\nu)$ is constraint to $n-d_i+\nu$ dimensions. Thus for the
last vertex $h(V_\delta)$ to intersect a line of given direction $\QQ
u_{(V_\delta,E_j)}$ and meeting a general $d_j$-codimensional affine
subspace $A_j$ requires $(n-d_i+\delta)+(n-d_j)\ge n$, that is,
$\delta\ge d_i+d_j-n$.
\qed
\medskip

Note that this observation can be refined by considering the
constraints imposed by the $A_i$ at any vertex.
\smallskip

By this argument and tedious case by case considerations we checked
the following examples.

\begin{example}
(1)\ \emph{Curves of degree $2$ in $\PP^3$ through $4$
points.}\ \ The expected number is zero because the dimension of the
space of homogeneous quadratic polynomials in two variables is $3$.
Hence every such curve lies in a plane, while this is not true for
$4$ general points.

The statement on the tropical side says that there are no tropical
curves of degree $2(1,0,0)$, $2(0,1,0)$, $2(0,0,1)$, $2(-1,-1,-1)$
through $4$ general points. The notation means that there are $2$
unbounded edges each in any of the directions $(1,0,0)$, $(0,1,0)$,
$(0,0,1)$, $(-1,-1,-1)$. Now there are very few marked trees
fulfilling the requirement of Proposition~\ref{distance of incidence
conditions}, and it turns out none of these could be the domain of a
tropical curve of the requested type. Alternatively, one can mimick
the argument on the geometric side by showing that each tropical
quadric lies in a tropical hyperplane, which does not contain $4$
general points.
\\[1ex]
(2)\ \emph{Curves of degree $(1,1,2)$ in $\PP^1\times\PP^1
\times\PP^1$ through $4$ points.}\ \ Again there is no such curve
because the projection to $\PP^1\times\PP^1$ is the graph of an
automorphism of $\PP^1$, which varies only in a three-dimensional
family.

Tropically we have degree $(\pm 1,0,0)$, $(0,\pm 1,0)$, $2(0,0,\pm 1)$,
so the projection to the first two coordinates is a tropical curve of
degree $(\pm 1,0)$, $(0,\pm 1)$ in $\QQ^2$. These form a
three-dimensional family, and hence there is no curve of this type
through $4$ general points.\\[1ex]
(3)\ \emph{Curves of degree $(1,2)$ in $\PP^1\times\PP^2$ through
$4$ points.}\ \ These are graphs of curves of degree $2$ in $\PP^2$,
and it is not hard to check that generically there is exactly one of
them. The corresponding computation on the affine side is for tropical
curves of degree $(0,0,\pm 1)$, $2(1,0,0)$, $2(0,1,0)$, $2(-1,-1,0)$.\\[1ex]
(4) \ \emph{Curves of degree $2$ in $\PP^3$ through $8$ lines, $5$ of
which non-general.}\ \  This example features higher dimensional
constraints, several types of tropical curves and $\tilde {\mathfrak
D}(\Gamma,\mathbf E, h,\mathbf A)\neq {\mathfrak D}(\Gamma,\mathbf E,
h,\mathbf A)$, that is, $\prod\delta_i\neq 1$. Consider quadrics in
$\PP^3$, so the degree is given by the vectors $(-1,0,0)$, $(0,-1,0)$,
$(0,0,-1)$, $(1,1,1)$, each with value $2$, so $e+n-3=8+3-3=8$. We
impose $8$ constraints by lines
\[\begin{array}{rclrcl}
L_1&=&(-2,-1,0)+\RR\cdot(1,\nu,0),&
L_2&=&(-1/2,-1/2,0)+\RR\cdot(\mu,-1,0),\\
L_3&=&(1/2,0,100)+\RR\cdot(-1,\lambda,1),\quad&
L_4&=&(-3,1,0)+\RR\cdot (0,0,1),\\
L_5&=&(-3,-1,0)+\RR\cdot (0,0,1),& L_6&=&(-1,-3,0)+\RR\cdot (0,0,1),\\
L_7&=&(1, -3,0)+\RR\cdot (0,0,1),&L_8&=&(3,2,0)+\RR\cdot (0,0,1),
\end{array}\]
with integers $\nu\ge 2$, $\mu\ge 1$, $\lambda\ge 3$, $\lambda\neq 6$.
Note the last five lines are parallel to each other. The composition
of a tropical curve $h:\Gamma\to\RR^3$ of given degree and matching
$\mathbf{A}$ with the projection $(x,y,z)\mapsto (x,y)$ onto the
$xy$-plane contracts the two unbounded edges in direction $(0,0,-1)$.
Hence it is a tropical plane quadric $\overline h:\overline\Gamma\to
\RR^2$ through the images of $L_4,\ldots L_8$, the $5$ points
$P_4=(-3,1)$, $P_5=(-3,-1)$, $P_6=(-1,-3)$, $P_7=(1,-3)$ and
$P_8=(3,2)$. There is exactly one such plane quadric, depicted in the
following figure.
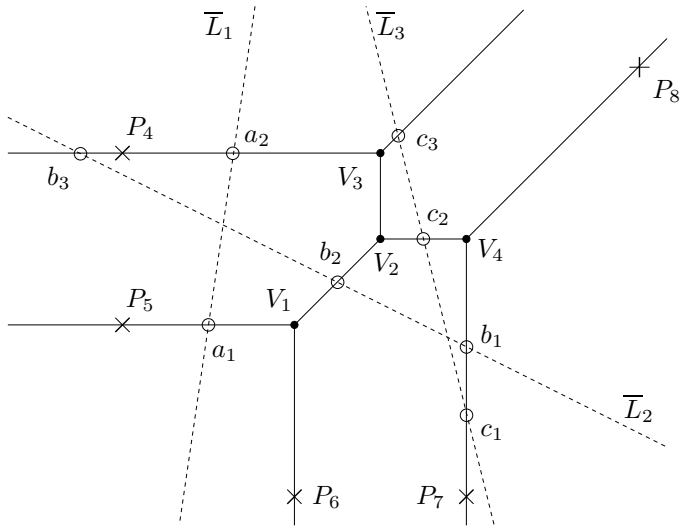
\begin{center}
\begin{picture}(0,0)%
\includegraphics{quadric2d.pstex}%
\end{picture}%
\setlength{\unitlength}{1184sp}%
\begingroup\makeatletter\ifx\SetFigFont\undefined%
\gdef\SetFigFont#1#2#3#4#5{%
  \reset@font\fontsize{#1}{#2pt}%
  \fontfamily{#3}\fontseries{#4}\fontshape{#5}%
  \selectfont}%
\fi\endgroup%
\begin{picture}(14845,10919)(-1821,-9983)
\put(676,-1786){\makebox(0,0)[lb]{\smash{{\SetFigFont{10}{12.0}{\familydefault}{\mddefault}{\updefault}{\color[rgb]{0,0,0}$P_4$}%
}}}}
\put(676,-5386){\makebox(0,0)[lb]{\smash{{\SetFigFont{10}{12.0}{\familydefault}{\mddefault}{\updefault}{\color[rgb]{0,0,0}$P_5$}%
}}}}
\put(4576,-9511){\makebox(0,0)[lb]{\smash{{\SetFigFont{10}{12.0}{\familydefault}{\mddefault}{\updefault}{\color[rgb]{0,0,0}$P_6$}%
}}}}
\put(6751,-9511){\makebox(0,0)[lb]{\smash{{\SetFigFont{10}{12.0}{\familydefault}{\mddefault}{\updefault}{\color[rgb]{0,0,0}$P_7$}%
}}}}
\put(3151,-1936){\makebox(0,0)[lb]{\smash{{\SetFigFont{10}{12.0}{\familydefault}{\mddefault}{\updefault}{\color[rgb]{0,0,0}$a_2$}%
}}}}
\put(2476,-6436){\makebox(0,0)[lb]{\smash{{\SetFigFont{10}{12.0}{\familydefault}{\mddefault}{\updefault}{\color[rgb]{0,0,0}$a_1$}%
}}}}
\put(-974,-2836){\makebox(0,0)[lb]{\smash{{\SetFigFont{10}{12.0}{\familydefault}{\mddefault}{\updefault}{\color[rgb]{0,0,0}$b_3$}%
}}}}
\put(8101,-6136){\makebox(0,0)[lb]{\smash{{\SetFigFont{10}{12.0}{\familydefault}{\mddefault}{\updefault}{\color[rgb]{0,0,0}$b_1$}%
}}}}
\put(8101,-8086){\makebox(0,0)[lb]{\smash{{\SetFigFont{10}{12.0}{\familydefault}{\mddefault}{\updefault}{\color[rgb]{0,0,0}$c_1$}%
}}}}
\put(6751,-2011){\makebox(0,0)[lb]{\smash{{\SetFigFont{10}{12.0}{\familydefault}{\mddefault}{\updefault}{\color[rgb]{0,0,0}$c_3$}%
}}}}
\put(5851,-4636){\makebox(0,0)[lb]{\smash{{\SetFigFont{10}{12.0}{\familydefault}{\mddefault}{\updefault}{\color[rgb]{0,0,0}$V_2$}%
}}}}
\put(8026,-4336){\makebox(0,0)[lb]{\smash{{\SetFigFont{10}{12.0}{\familydefault}{\mddefault}{\updefault}{\color[rgb]{0,0,0}$V_4$}%
}}}}
\put(5101,-2836){\makebox(0,0)[lb]{\smash{{\SetFigFont{10}{12.0}{\familydefault}{\mddefault}{\updefault}{\color[rgb]{0,0,0}$V_3$}%
}}}}
\put(3601,-5386){\makebox(0,0)[lb]{\smash{{\SetFigFont{10}{12.0}{\familydefault}{\mddefault}{\updefault}{\color[rgb]{0,0,0}$V_1$}%
}}}}
\put(4726,-4411){\makebox(0,0)[lb]{\smash{{\SetFigFont{10}{12.0}{\familydefault}{\mddefault}{\updefault}{\color[rgb]{0,0,0}$b_2$}%
}}}}
\put(6976,-3586){\makebox(0,0)[lb]{\smash{{\SetFigFont{10}{12.0}{\familydefault}{\mddefault}{\updefault}{\color[rgb]{0,0,0}$c_2$}%
}}}}
\put(2326,314){\makebox(0,0)[lb]{\smash{{\SetFigFont{10}{12.0}{\familydefault}{\mddefault}{\updefault}{\color[rgb]{0,0,0}$\overline L_1$}%
}}}}
\put(5926,314){\makebox(0,0)[lb]{\smash{{\SetFigFont{10}{12.0}{\familydefault}{\mddefault}{\updefault}{\color[rgb]{0,0,0}$\overline L_3$}%
}}}}
\put(11101,-7711){\makebox(0,0)[lb]{\smash{{\SetFigFont{10}{12.0}{\familydefault}{\mddefault}{\updefault}{\color[rgb]{0,0,0}$\overline L_2$}%
}}}}
\put(11701,-1017){\makebox(0,0)[lb]{\smash{{\SetFigFont{10}{12.0}{\familydefault}{\mddefault}{\updefault}{\color[rgb]{0,0,0}$P_8$}%
}}}}
\end{picture}%

\end{center}
\end{example}
The images of $L_1$, $L_2$, $L_3$ are lines $\overline L_1$,
$\overline L_2$, $\overline L_3$. The possible images of the first
three marked points are denoted $a_i$, $b_j$, $c_k$. It is not hard
to see that each of the $2\cdot 3\cdot3=18$ possibilities $(a_i b_j
c_k)$ gives rise to a unique lift of $\overline h$ to a tropical
quadric in $\RR^3$. Because on the considered region $L_3$ has a
large value of the last coordinate the two contracted unbounded edges
always lie over the selected $a_i$ and $b_j$, and these unbounded
edges are the first and second marked edges.

It remains to compute the contribution to $N_{0,\Delta}^\text{trop}
(\mathbf A)$ for each of the $18$ tropical curves thus obtained.  We
explain the calculation in the case ($a_1 b_3 c_j$), $j = 1, 2, 3$.
In this case, the vertices of the tropical curve are $V_1, V_2, V_3,
V_4$ and lifts $A_1, B_3$ of $a_1$, $b_3$. ${\mathfrak D}(\Gamma,\mathbf E, h,
\mathbf A)$ is the index of the inclusion $\Phi$ of lattices
defined in (\ref{lattice map}) (Proposition~\ref{count of maximally
degenerate curves}). The domain of $\Phi$ is a copy of $\ZZ^3$ for
each of the $6$~vertices. Define primitive vectors $f_i\in\ZZ^3$ in
the direction of the images of the bounded edges as in the figure
below. Now the range of $\Phi$ is the free part of the following
abelian group.
\[
\begin{array}{l}
\ZZ^3 /\ZZ f_1 \oplus 
\ZZ^3 /\ZZ  f_2 \oplus 
\ZZ^3 /\ZZ  f_3 \oplus 
\ZZ^3 /\ZZ  f_4 \oplus 
\ZZ^3 /\ZZ  f_5 \\
\oplus\, \ZZ^3 / (\ZZ  (e_1+\nu e_2) +\ZZ  {e_3}) \oplus 
\ZZ^3 / (\ZZ  (\mu e_1 - e_2) +\ZZ  {e_3}) \oplus 
\ZZ^3 / (\ZZ  (-e_1 + \lambda e_2 + e_3) +\ZZ  f)\\
\oplus\, \ZZ^3 / \langle e_1,e_3\rangle \oplus
\ZZ^3 / \langle e_1,e_3\rangle \oplus
\ZZ^3 / \langle e_2,e_3\rangle \oplus
\ZZ^3 / \langle e_2,e_3\rangle \oplus
\ZZ^3 / \langle (1,1,1),e_3\rangle.
\end{array}
\]
Here ${e_1} = (1, 0, 0)$, ${e_2} = (0, 1, 0)$, ${e_3} = (0, 0, 1)$.
The vector $f$ in the last summand is ${(0, 1, 0)}$ for $(a_1 b_3
c_1)$, ${(1, 0, 1)}$ for $(a_1 b_3 c_2)$, and  ${(1, 1, 1)}$ for
$(a_1 b_3 c_3)$. This holds regardless of the positions of $b_1$ and
$c_1$ relative to $P_7$ and of $b_3$ relative to $P_4$, which depend
on the choices of $\lambda$ and $\mu$.  The image of $(z_1, \cdots,
z_6) \in (\ZZ^3)^{\oplus 6}$ under $\Phi$ is
\begin{equation*}
\begin{array}{l}
(z_2 - z_1, z_3 - z_2, z_4 - z_2, z_1 - z_5, z_3 - z_6, \\
z_5 - r_1, z_6 - r_2, z - r_3,\\
z_6 - r_4, z_5 - r_5, z_1 - r_6, z_4 - r_7, z_4 - r_8),
\end{array}
\end{equation*}
where $r_i\in L_i\cap \ZZ^3$ are arbitrary and each term is to be
viewed in the quotient group. $z$ in the last term is given by $z =
z_4$ for $(a_1 b_3 c_1)$ or $(a_1 b_3 c_2)$, and $z = z_3$ for $(a_1
b_3 c_3)$. The contribution $\tilde {\mathfrak D}(\Gamma,\mathbf
E,h,\mathbf A)$ is given by the absolute value of the determinant of
the $18 \times 18$-matrix representing the linear part of this map,
times the order of the torsion group of the abelian group above. This
gives $\nu$, $\lambda \nu$ and $(1+\lambda)\nu$ for $(a_1 b_3 c_1)$,
$(a_1 b_3 c_2)$, $(a_1 b_3 c_3)$, respectively.
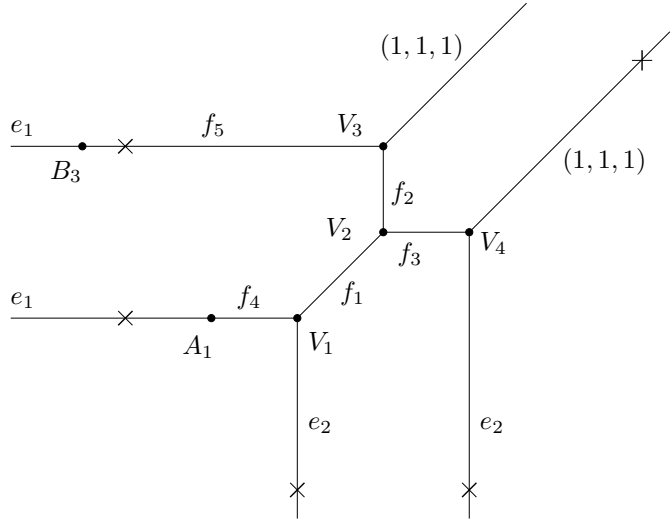
\begin{center}
\begin{picture}(0,0)%
\includegraphics{quadric3d.pstex}%
\end{picture}%
\setlength{\unitlength}{1184sp}%
\begingroup\makeatletter\ifx\SetFigFont\undefined%
\gdef\SetFigFont#1#2#3#4#5{%
  \reset@font\fontsize{#1}{#2pt}%
  \fontfamily{#3}\fontseries{#4}\fontshape{#5}%
  \selectfont}%
\fi\endgroup%
\begin{picture}(13844,10844)(-1821,-9983)
\put(8026,-4336){\makebox(0,0)[lb]{\smash{{\SetFigFont{10}{12.0}{\familydefault}{\mddefault}{\updefault}{\color[rgb]{0,0,0}$V_4$}%
}}}}
\put(4426,-6436){\makebox(0,0)[lb]{\smash{{\SetFigFont{10}{12.0}{\familydefault}{\mddefault}{\updefault}{\color[rgb]{0,0,0}$V_1$}%
}}}}
\put(-974,-2836){\makebox(0,0)[lb]{\smash{{\SetFigFont{10}{12.0}{\familydefault}{\mddefault}{\updefault}{\color[rgb]{0,0,0}$B_3$}%
}}}}
\put(5926,-211){\makebox(0,0)[lb]{\smash{{\SetFigFont{10}{12.0}{\familydefault}{\mddefault}{\updefault}{\color[rgb]{0,0,0}$(1,1,1)$}%
}}}}
\put(5026,-1861){\makebox(0,0)[lb]{\smash{{\SetFigFont{10}{12.0}{\familydefault}{\mddefault}{\updefault}{\color[rgb]{0,0,0}$V_3$}%
}}}}
\put(6301,-4561){\makebox(0,0)[lb]{\smash{{\SetFigFont{10}{12.0}{\familydefault}{\mddefault}{\updefault}{\color[rgb]{0,0,0}$f_3$}%
}}}}
\put(-1799,-5461){\makebox(0,0)[lb]{\smash{{\SetFigFont{10}{12.0}{\familydefault}{\mddefault}{\updefault}{\color[rgb]{0,0,0}$e_1$}%
}}}}
\put(-1799,-1861){\makebox(0,0)[lb]{\smash{{\SetFigFont{10}{12.0}{\familydefault}{\mddefault}{\updefault}{\color[rgb]{0,0,0}$e_1$}%
}}}}
\put(1801,-6511){\makebox(0,0)[lb]{\smash{{\SetFigFont{10}{12.0}{\familydefault}{\mddefault}{\updefault}{\color[rgb]{0,0,0}$A_1$}%
}}}}
\put(2926,-5461){\makebox(0,0)[lb]{\smash{{\SetFigFont{10}{12.0}{\familydefault}{\mddefault}{\updefault}{\color[rgb]{0,0,0}$f_4$}%
}}}}
\put(2176,-1861){\makebox(0,0)[lb]{\smash{{\SetFigFont{10}{12.0}{\familydefault}{\mddefault}{\updefault}{\color[rgb]{0,0,0}$f_5$}%
}}}}
\put(5101,-5386){\makebox(0,0)[lb]{\smash{{\SetFigFont{10}{12.0}{\familydefault}{\mddefault}{\updefault}{\color[rgb]{0,0,0}$f_1$}%
}}}}
\put(6151,-3286){\makebox(0,0)[lb]{\smash{{\SetFigFont{10}{12.0}{\familydefault}{\mddefault}{\updefault}{\color[rgb]{0,0,0}$f_2$}%
}}}}
\put(4801,-4036){\makebox(0,0)[lb]{\smash{{\SetFigFont{10}{12.0}{\familydefault}{\mddefault}{\updefault}{\color[rgb]{0,0,0}$V_2$}%
}}}}
\put(4426,-8086){\makebox(0,0)[lb]{\smash{{\SetFigFont{10}{12.0}{\familydefault}{\mddefault}{\updefault}{\color[rgb]{0,0,0}$e_2$}%
}}}}
\put(8026,-8086){\makebox(0,0)[lb]{\smash{{\SetFigFont{10}{12.0}{\familydefault}{\mddefault}{\updefault}{\color[rgb]{0,0,0}$e_2$}%
}}}}
\put(9751,-2611){\makebox(0,0)[lb]{\smash{{\SetFigFont{10}{12.0}{\familydefault}{\mddefault}{\updefault}{\color[rgb]{0,0,0}$(1,1,1)$}%
}}}}
\end{picture}%

\end{center}
Doing the other cases $(a_i b_j c_k)$ it turns out the result never
depends on the choice of $a_i$. The following table now lists
$\tilde {\mathfrak D}(\Gamma,\mathbf E,h,\mathbf A)$ for all $18$ cases.
\begin{center}
{\small
\begin{tabular}{| l | l | l | l | l | l | l | l | l |}
\hline
$a_i b_1 c_1$ & $a_i b_1 c_2$ & $a_i b_1 c_3$ & $a_i b_2 c_1$ &
$a_i b_2 c_2$ &  $a_i b_2 c_3$ & $a_i b_3 c_1$ & $a_i b_3 c_2$ &
$a_i b_3 c_3$ \\ \hline
$\mu\nu$ & $\lambda \mu\nu$ & $(1+\lambda)\mu\nu$
& $(1+\mu)\nu$ & $\lambda(1+\mu)\nu$ & $(1+\lambda)(1+\mu)\nu$
& $\nu$ & $\lambda \nu$
& $(1+\lambda)\nu$\\ \hline
\end{tabular}}
\end{center}
If $\lambda$ is odd $\tilde {\mathfrak D}(\Gamma,\mathbf E,h,\mathbf
A)= {\mathfrak D}(\Gamma,\mathbf E,h,\mathbf A)$ for all $(a_i b_j
c_k)$; but if $\lambda$ is even the last summand in the definition of
the range of $\Phi$ for the $(a_i b_3 c_2)$-case, that is, $\ZZ^3 /
(\ZZ(-e_1 - \lambda e_2 + e_3)+ \ZZ f)$, $f = (1, 0, 1)$, is
isomorphic to $\ZZ \oplus \ZZ/2\ZZ$. In fact, in this case ${\mathfrak
D}(\Gamma,\mathbf E,h, \mathbf A)=\lambda\nu/2$, $\prod_i\delta_i=2$
and $\tilde {\mathfrak D}(\Gamma,\mathbf E,h, \mathbf A)=\lambda\nu$.

Finally, taking twice (for $i=1,2$) the sum of
the terms in the table gives
\[
N_{0,\Delta}^\text{trop}(\mathbf A)= 8\nu(\mu+1)(\lambda+1).
\]

This number can be checked by a direct algebraic-geometric
computation, say for $\kk=\CC$. The five parallel lines define orbits
of the form $(\alpha,\beta,\gamma t)\in (\CC^*)^3$, $t\in \CC^*$.
Their closures are lines through $P=[0,0,1,0]\in\PP^3$. Projecting
from this point defines a map $\pi:\PP^3\setminus \{P\}\to \PP^2$
contracting these lines to five points. For a general choice of lines
through $P$ there is a unique quadric $Q\subset \PP^2$ passing
through them, which is smooth. The closure $\tilde Q$ of
$\pi^{-1}(Q)$ is then a quadric cone with isolated singularity at
$P$. The intersection of $\tilde Q$ with a hyperplane $H\subset
\PP^3$ defines a quadric curve intersecting the five given lines.
Conversely, any quadric curve in $\PP^3$ lies in a hyperplane, as
already mentioned in the first example. Thus \emph{all} quadric
curves intersecting the five lines in $(\CC^*)^3\subset \PP^3$ arise
in this way.

Now the remaining three incidence curves $Z_1,Z_2,Z_3$ are rational
curves of degrees $\nu$, $\mu+1$ and $\lambda+1$. Hence they
intersect $\tilde Q$ in $2\nu$, $2(\mu+1)$ and $2(\lambda+1)$ points,
respectively, and all these points are disjoint for a general choice
of $Z_1,Z_2,Z_3$. Each of the $2\nu\cdot 2(\mu+1)\cdot 2(\lambda+1)$
choices of three intersection points determines a unique hyperplane
$H\subset \PP^3$. Then $H\cap\tilde Q$ is the requested quadric
curve intersecting $Z_1,\ldots Z_8$. Thus
\[
N_{0,\Delta}^\text{alg} (\mathbf L)= 8\nu(\mu+1)(\lambda+1),
\]
with $\mathbf L=(L(A_1),\ldots,L(A_8))$, in agreement with the tropical
computation.
\qed


\end{document}